\documentclass[EJP]{ejpecp} 

\usepackage[utf8]{inputenc}

\allowdisplaybreaks

\SHORTTITLE{Particle representations for SPDEs with boundary conditions} 

\TITLE{Particle representations for stochastic partial differential equations with boundary conditions} 

\AUTHORS{%
Dan~Crisan \footnote{Department of Mathematics, Imperial College London. 180 Queen's Gate, London SW7 2AZ, United Kingdom. \EMAIL{d.crisan@imperial.ac.uk}. \url{http://www.imperial.ac.uk/\~dcrisan/}.  Research partially supported by an EPSRC Mathematics Platform grant.}
  \and 
  Christopher~Janjigian\footnote{Department of Mathematics, University of Utah. 155 South 1400 East, Salt Lake City, UT  84112. \EMAIL{janjigia@math.utah.edu}. \url{http://www.math.utah.edu/\~janjigia/}. Research supported in part by NSF grant DMS 11-06424 and in part by a postdoctoral grant from the Fondation Sciences Math\'ematiques de Paris.}
  \and 
  Thomas~G.~Kurtz\footnote{Departments of Mathematics and Statistics, University of Wisconsin - Madison. 480 Lincoln Drive, Madison, WI  53706-1388. \EMAIL{kurtz@math.wisc.edu}. \url{http://www.math.wisc.edu/\~kurtz/}. Research supported in part by NSF grant
DMS 11-06424 and by a Nelder Fellowship at Imperial College London.}
}


\AMSSUBJ{Primary:  
60H15; 60H35 Secondary: 60B12; 60F17; 60F25;
60H10; 35R60; 93E11} 

\SUBMITTED{July 29, 2016} 
\ACCEPTED{June 7, 2018} 

\ARXIVID{1607.08909} 

\VOLUME{0}
\YEAR{2016}
\PAPERNUM{0}
\DOI{10.1214/YY-TN}

\ABSTRACT{In this article, we study weighted particle 
representations for a class of stochastic partial 
differential equations (SPDE) with Dirichlet boundary 
conditions.  The locations and weights of the particles 
satisfy an infinite system of stochastic differential 
equations.  The locations are given by independent, 
stationary reflecting diffusions in a bounded domain, and 
the weights evolve according to an infinite system of 
stochastic differential equations driven by a 
\emph{common} Gaussian white noise $W$ which is the 
stochastic input for the SPDE.  The weights interact 
through $V$, the associated weighted empirical measure, 
which gives the solution of the SPDE.  When a particle 
hits the boundary its weight jumps to a value given by 
a function of the location of the particle on the 
boundary.  This function determines the boundary 
condition for the SPDE.  We show existence and 
uniqueness of a solution of the infinite system of 
stochastic differential equations giving the locations and 
weights of the particles and derive two weak forms for 
the corresponding SPDE depending on the choice of test 
functions.  The weighted empirical measure $V$ is the 
unique solution for each of the nonlinear stochastic 
partial differential equations.  The work is motivated by 
and applied to the stochastic Allen-Cahn equation and 
extends the earlier of work of Kurtz and Xiong in 
\cite{KX99,KX04}.  }

\renewcommand {\theequation}{\arabic{section}.\arabic{equation}}
\def\non{\nonumber}

\def\bar{\overline}

\begin{document}



\setcounter{equation}{0}

\section{Introduction} In the following, we study particle 
representations for a class of nonlinear stochastic 
partial differential equations that includes the stochastic 
version of the Allen-Cahn equation \cite{ac72,ac73} and 
that of the equation governing the stochastic 
quantization of $\Phi^4_d$ Euclidean quantum field theory with 
quartic interaction \cite{pw}, that is, the equation 
\begin{equation}\label{aceq0}dv=\Delta v+G(v)v+W,\end{equation}
where $G$ is a (possibly) nonlinear function
\footnote{In the 
original work of Allen and Cahn (see \cite{ac72,ac73}), 
$G(v)=1-v^2$ for all $v\in\mathbb{R}$ whilst in the case of 
the $\Phi^4_d$ equation of Euclidean quantum field theory, 
$G(v)=-v^2$ ($d$ represents the state space dimension).} 
and 
$W$ is a space-time noise
\footnote{A detailed description 
of the noise $W$ is given below and in the Appendix.}.  
These particle representations lead naturally to the 
solution of a weak version of a stochastic partial 
differential equation similar to (\ref{aceq0}).  See 
equation (\ref{aceq}) below and Section \ref{secteqbc}.  

The  approach taken here has its roots in the study of 
the McKean-Vlasov problem and its stochastic 
perturbation.  In its simplest form, the problem begins 
with a finite system of stochastic differential equations
\begin{equation}X_i^n(t)=X_i^n(0)+\int_0^t\sigma (X_i^n(s),V^n(s)
)dB_i(s)+\int_0^tc(X_i^n(s),V^n(s))ds,\quad 1\leq i\leq n,\label{fisys}\end{equation}
where the $X_i^n$ take values in ${\Bbb R}^d$, $V^n(t)$ is the empirical 
measure $\frac 1n\sum_{i=1}^n\delta_{X_i^n(t)}$, and the $B_i$ are independent, 
standard Brownian motions in an appropriate Euclidean 
space.  The primary goal in this setting 
is to prove that the sequence of 
empirical measures $V^n$ converges in distribution and to 
characterize the limit $V$ as a measure valued process 
which solves the following nonlinear partial differential 
equation, written in weak form, \footnote{In 
(\ref{mkpde}) and thereafter, we use the notation $\langle\varphi 
,\mu\rangle$ 
to express the integral of $\varphi$ with respect to $\mu$, that is, 
$\langle\varphi ,\mu\rangle =\int\varphi (x)\mu (dx).$} 
\begin{equation}\langle\varphi ,V(t)\rangle =\langle\varphi ,V(0)
\rangle +\int_0^t\langle L(V(s))\varphi ,V(s)\rangle ds,\label{mkpde}\end{equation}
where $\varphi :{\Bbb R}^d\mapsto {\Bbb R}$ belongs to a suitably chosen class of Borel measurable functions and 
\[\langle\varphi ,V(t)\rangle =\int\varphi dV(t)=\int\varphi (u)V
(t,du).\]
In (\ref{mkpde}), we use $a(x,\nu )=\sigma (x,\nu )\sigma (x,\nu )^T$ and  $L(\nu)$ is the 
differential operator
\[L(\nu )\varphi (x)=\frac 12\sum_{i,j}a_{ij}(x,\nu )\partial^2_{
x_i}{}_{x_j}\varphi (x)+\sum_ic_i(x,\nu )\partial_{x_i}\varphi (x
).\]
There are many approaches to this problem
 \cite{Gra92,McK67,Oel84}. (See also the recent book 
\cite{Kol10}.)  The approach in which we are interested, 
introduced in \cite{KP96b} and developed further in 
\cite{CKL14,KK10,KX99}, is simply to let the limit be 
given by the infinite system
\begin{equation}X_i(t)=X_i(0)+\int_0^t\sigma (X_i(s),V(s))dB_i(s)
+\int_0^tc(X_i(s),V(s))ds,\quad 1\leq i<\infty .\label{infsys}\end{equation}
To make sense out of this system (in particular, the 
relationship of $V$ to the $X_i$), 
note that we can 
assume, without loss of generality, that the finite 
system $\{X_i^n\}$ is exchangeable (randomly permute the index 
$i$), so if one shows relative compactness of the sequence, 
any limit point will be an infinite exchangeable 
sequence and we can require $V(t)$ to be the de Finetti 
measure for the sequence $\{X_i(t)\}$, that is,
\[V(t)=\lim_{m\rightarrow\infty}\frac 1m\sum_{i=1}^m\delta_{X_i(t
)},\]
in the sense that $\langle\varphi ,V(t)\rangle =\lim_{m\rightarrow
\infty}\frac 1m\sum_{i=1}^m\varphi (X_i(t))$ for all 
bounded, measurable $\varphi$.  More precisely, $V$ is a process 
with sample paths in $C_{{\cal P}({\Bbb R}^d)}[0,\infty )$. (See Lemma 4.4 of 
\cite{KK10}.)

Note that while the $X_i^n$ give a {\em particle approximation }
of the solution of (\ref{mkpde}), the $X_i$ give a {\em particle }
{\em representation\/} of the solution, that is, the de Finetti 
measure of $\{X_i(t)\}$ is the desired $V(t)$.

In the following, we will make use of representations 
similar to those studied in \cite{KX99}.  They differ 
from those considered in the other papers mentioned as 
each particle is given a weight $A_i(t)$ and $V$ is given by 
\[V(t)=\lim_{m\rightarrow\infty}\frac 1m\sum_{i=1}^mA_i(t)\delta_{
X_i(t)}.\]

The models in the current paper differ from those in 
\cite{KX99} in two primary ways.  First, the $X_i$ will be 
independent, stationary diffusion processes defined on a 
domain $D\subset\mathbb{R}^d$ with reflecting boundary.  The 
stationary distribution will be denoted by $\pi$.  Second, we 
will place boundary conditions on the solution.  
Essentially, if $v(t,u)$ is an appropriately defined solution 
to a stochastic partial differential equation of interest, 
we would like to have 
\begin{equation}v(t,u)=g(u),\quad u\in\partial D,\;t>0.\label{dirchbc}\end{equation}
The precise sense in which the boundary 
conditions hold will be discussed later and
will depend on the conditions assumed.

Specifically, in the same vein as equation (\ref{mkpde}), we will 
consider a class of nonlinear stochastic partial 
differential equations  
written in weak form 
\begin{eqnarray}
\langle\varphi ,V(t)\rangle =\langle\varphi ,V(0)\rangle +\int_0^
t\langle {\Bbb L}\varphi ,V(s)\rangle ds+\int_0^t\langle G(v(s,\cdot 
),\cdot )\varphi ,V(s)\rangle ds\nonumber\\
+\int_0^t\int_{\bar {D}}\varphi (x)b(x)\pi (dx)ds+\int_{{\Bbb U}\times 
[0,t]}\int_{\bar {D}}\varphi (x)\rho (x,u)\pi (dx)W(du,ds),
\label{aceq}\end{eqnarray}
where $D$ is a bounded, open, connected subset of ${\Bbb R}^d$,
for the moment  the test functions are 
$\varphi\in C^2_c(D)$ 
the twice continuously differentiable 
functions with compact support in  $D$, 
\[{\Bbb L}\varphi (x)=\frac 12\sum_{i,j}a_{ij}(x)\partial^2_{x_ix_
j}\varphi (x)+\sum_ic_i(x)\partial_{x_i}\varphi (x),\]
$v(s,x)$ is the density of $V(s)$ with respect to $\pi$, and 
 $\langle\varphi ,V(0)\rangle =\int\varphi (x)h(x)\pi (dx)$ for a specified $
h$. Note that if 
the solution $V$ is adapted to a filtration $\{{\cal F}_t\}$, then 
we can assume that $v$ is a progressively measurable 
process with values in $L^1(\pi )$. In particular, the mapping 
$(s,x,\omega )\rightarrow v(s,x,\omega )$ is ${\cal B}[0,\infty )
\times {\cal B}(\bar {D})\times\vee_t{\cal F}_t$ measurable. (${\cal B}
(E)$ 
denotes the Borel subsets of a metric space $E$.)  
The existence of a measurable version of the
 process of densities follows by a monotone class argument. See 
Lemma \ref{opmeas} for a similar argument.

Throughout, we will assume that  ${\Bbb U}$ is a complete, 
separable metric space, $\mu$ is a $\sigma$-finite Borel measure on 
${\Bbb U}$, and $\ell$ is Lebesgue measure on $[0,\infty )$.
$W$ is Gaussian white noise on ${\Bbb U}\times [0,\infty )$ with covariance 
measure $\mu\times\ell$, that is, $W(C\times [0,t])$ has expectation zero 
for all $0\leq t<\infty$ and $C\in {\cal B}({\Bbb U})$ with $\mu 
(C)<\infty$, and 
${\Bbb E}[W(C_1\times [0,t])W($$C_2\times [0,s])$$=t\wedge s\mu (
C_1\cap C_2)$.
(See Appendix \ref{bcapp1}.)

Formally, equation (\ref{aceq}) is the weak form of 
\begin{equation}\label{sf}v(t,x)=v(0,x)+\int_0^t\left[\mathbb L^{*}v(s,
x)+v(s,x)G(v(s,x),x)+b(x)\right]ds+\int_{{\Bbb U}\times [0,t]}\rho (x,
u)W(du,ds),\end{equation}
where $\mathbb L^{*}$ is the formal adjoint of the operator 
$\mathbb L$ with respect to the inner product $\langle\psi ,\varphi
\rangle_{\pi}=\int_{\bar {D}}\psi (x)\varphi (x)\pi (dx)$.

To obtain, for example, the stochastic 
Allen-Cahn equation (\ref{aceq0}), we can choose ${\Bbb L}=\Delta$.

\vspace{.1in}
  
\setcounter{equation}{0} \section{Basic conditions and 
statement of main theorems} With the results of 
\cite{KX99} in mind, we are interested in solutions $V$ of 
(\ref{aceq}) that are measures or perhaps signed 
measures.  We emphasize here that we are talking about 
a {\em representation\/} of the solution of the equation{\em ,\/} not a 
limit or approximation theorem (although these 
representations can be used to prove limit theorems).  
To obtain the representation of $V(t)$ we desire, we must 
identify both the locations $X_i(t)$ and the weights $A_i(t)$.  
The sequence $\{(X_i(t),A_i(t))\}$ is required to be 
exchangeable, so by de Finetti's theorem, for each $t$ 
there will exist a random measure $\Xi (t)$, which we will 
refer to as the de Finetti measure for the sequence 
$\{(X_i(t),A_i(t))\}$, such that for each bounded, measurable 
function $\psi$ on $\bar {D}\times {\Bbb R}$, 
\[\langle\psi ,\Xi (t)\rangle =\lim_{m\rightarrow\infty}\frac 1m\sum_{
i=1}^m\psi (X_i(t),A_i(t))=\int_{\bar {D}\times {\Bbb R}}\psi (x,
a)\Xi (t,dx,da).\]
Then, assuming ${\Bbb E}[|A_i(t)|]<\infty$, $V(t)$ will be given by
\begin{equation}\langle\varphi ,V(t)\rangle =\lim_{m\rightarrow\infty}\frac 
1m\sum_{i=1}^mA_i(t)\varphi (X_i(t))=\int_{\bar {D}\times {\Bbb R}}
a\varphi (x)\Xi (t,dx,da).\label{solrep}\end{equation}
If the $A_i(t)$ are nonnegative, then $V(t)$ will be a measure, 
but we do not rule out the possibility that the $A_i(t)$ can be 
negative and $V(t)$ a signed measure.
The weights and locations will be solutions of an infinite 
system of stochastic differential equations that are 
coupled only through $V$ and common noise terms.

We take the $X_i$ to be independent, stationary solutions 
of the
Skorohod equation
\begin{equation}X_i(t)=X_i(0)+\int_0^t\sigma (X_i(s))dB_i(s)+\int_
0^tc(Xi(s))ds+\int_0^t\eta (X_i(s))dL_i(s),\ i\ge 1,\label{infsys2}\end{equation}
where $\eta (x)$ is a vector field defined on 
the boundary 
$\partial D$, and $L_i$ is a local time on $\partial D$ for $X_i$, that is, $
L_i$ is
a nondecreasing process that increases 
only when $X_i$ is in $\partial D$. To avoid some of the 
complexities of reflecting diffusions and focus on the 
new ideas in our representation, we will assume the 
following condition throughout.

\begin{condition}\label{diffcnd}
\begin{itemize}
\item[a)] $D\subset {\Bbb R}^d$ is bounded, open,  connected, and has 
a $C^2$ boundary.  

\item[b)] $\sigma$, $c$, and $\eta$ are continuous, $\sigma$ is nondegenerate 
on $\bar {D}$, and $\eta (x)\cdot n_D(x)>0$, $x\in\partial D$, where $
n_D(x)$ is the unit
inward normal at $x$.  

\item[c)] For a standard Brownian motion $B$ and $X(0)\in\bar {D}$ 
independent of $B$, the solution 
of
\begin{equation}X(t)=X(0)+\int_0^t\sigma (X(s))dB(s)+\int_0^tc(X(
s))ds+\int_0^t\eta (X(s))dL(s),\ \label{rsde}\end{equation}
 is (weakly) unique for all initial 
distributions in ${\cal P}(\bar {D})$. 

\item[d)] The $X_i$ are independent
solutions of (\ref{infsys2}) with independent standard 
Brownian motions $B_i$ 
and independent and identically distributed  $X_i(0)\in\bar {D}$. 
The distribution $\pi$ of $X_i(0)$ is a stationary distribution 
for (\ref{rsde}).
\end{itemize}
\end{condition}

\begin{remark}
See \cite{DI93} for conditions implying strong uniqueness 
for (\ref{rsde}).  Weak uniqueness can be proved 
employing results from partial differential 
equations. (See, for example, Theorem 8.1.5 in \cite{EK86}. 
In the notation of that 
theorem, $\eta (x)=-c(x)$.) 
Weak uniqueness can also be obtained via 
submartingale problems. (See \cite{SV71}.)  These 
approaches all require different conditions on $D$ and the 
coefficients, so we simply assume uniqueness.
\end{remark}

\begin{lemma}\label{stcsg}
The equation (\ref{rsde}) determines a strong Markov 
process, and 
the corresponding semigroup defined by
\[{\Bbb T}(t)\varphi (x)={\Bbb E}[\varphi (X(t))|X(0)=x]\]
satisfies 
\begin{equation}{\Bbb T}(t):C_b(\bar {D})\rightarrow C_b(\bar {D}
).\label{fellprop}\end{equation}
\end{lemma}

\begin{proof}
Let $X^x$ denote the solution of 
(\ref{rsde}) with $X^x(0)=x$.  Continuity of 
the coefficients and the assumption of uniqueness implies 
that the mapping 
$x\rightarrow X^x$ is continuous in the sense of convergence in 
distribution in $C_{\bar {D}}[0,\infty )$. The simplest way to see this 
continuity is 
to first define the time change 
\[\tau^x(t)=\inf\{r;r+L^x(r)>t\}.\]
Then setting $\lambda^x(t)=L^x(\tau^x(t))$ and $Y^x(t)=X^x(\tau^x
(t))$, we have 
\[Y^x(t)=x+\int_0^t\sigma (Y^x(s))dB\circ\tau^x(s)+\int_0^tc(Y^x(
s))ds+\int_0^t\eta (Y^x(s))d\lambda^x(s).\]

Since all the coefficients are continuous and 
$\tau^x(t)+\lambda^x(t)=t$, $\tau^x$ and $\lambda^x$ are Lipschitz with constant $
1$, 
$\{Y^x,x\in\bar {D}\}$ is relatively compact, and any limit point 
$(Y^{x_0},\tau^{x_0},\lambda^{x_0})$ of $(Y^x,\tau^x,\lambda^x)$ with $
x\rightarrow x_0$ will satisfy 
\[Y^{x_0}(t)=x_0+\int_0^t\sigma (Y^{x_0}(s))dB\circ\tau^{x_0}(s)+
\int_0^tc(Y^{x_0}(s))ds+\int_0^t\eta (Y^{x_0}(s))d\lambda^{x_0}(s
).\]
Since $\lambda^{x_0}$ increases only when $Y^{x_0}\in\partial D$, the assumptions 
on $\eta$ ensure that $\tau^{x_0}$ is strictly increasing (as are the 
$\tau^x$).  Letting $\gamma^{x_0}$ denote the inverse of $\tau^{x_
0}$, 
$X^{x_0}\equiv Y^{x_0}\circ\gamma^{x_0}$, is a solution of (\ref{rsde}) with 
$X(0)=x_0$, and by uniqueness, $X^x\Rightarrow X^{x_0}$ as 
$x\rightarrow x_0$ giving (\ref{fellprop}).
\end{proof}

  If $\varphi\in C_b^2(\bar {D})$, then by It\^o's formula,
\begin{eqnarray}
\varphi (X(t))=\varphi (X(0))+\int_0^t\nabla\varphi (X(s))^T\sigma 
(X(s))dB(s)+\int_0^t{\Bbb L}\varphi (X(s))ds\label{itophi}\\
+\int_0^t\nabla\varphi (X(s))\cdot\eta (X(s))dL(s),\ \nonumber\end{eqnarray}
where
\begin{equation}{\Bbb L}\varphi (x)=\frac 12\sum_{i,j}a_{ij}(x)\partial^2_{x_ix_
j}\varphi (x)+\sum_ic_i(x)\partial_{x_i}\varphi (x),\label{infgen}\end{equation}
with  $a(x)=\sigma (x)\sigma (x)^T$,
where $\sigma^T$ is the transpose of $\sigma$.

\begin{lemma}
The infinitesimal generator ${\Bbb A}$ for the 
semigroup $\{{\Bbb T}(t)\}$ is an extension of
\begin{equation}\{(\varphi ,{\Bbb L}\varphi ):\varphi\in C^2_b(\bar {
D}),\;\eta (x)\cdot\nabla\varphi |_{\partial D}=0\}.\label{gendef}\end{equation}
\end{lemma}

\begin{proof}
The boundary condition and the martingale property of 
the stochastic integral implies
\begin{equation}\frac {{\Bbb T}(t)\varphi (x)-\varphi (x)}t=\frac 
1t\int_0^t{\Bbb E}[{\Bbb L}\varphi (X^x(s))]ds.\label{sggen}\end{equation}
The continuity in distribution of $x\rightarrow X^x$ and the 
continuity of ${\Bbb L}\varphi$ assure that if $x\rightarrow x_0$ and $
t\rightarrow 0$, the right 
side of (\ref{sggen}) converges to ${\Bbb L}\varphi (x_0)$.  The 
compactness of $\bar {D}$ then assures that the convergence of 
the left side of (\ref{sggen}) to ${\Bbb L}\varphi (x)$ is uniform in $
x$.
\end{proof}

\begin{lemma}
Let $X,L$ be a solution of (\ref{rsde}).  Then for each 
$t>0$,
\[\int_0^t{\bf 1}_{\partial D}(X(s))ds=0\quad a.s.,\quad t\geq 0,\]
and ${\Bbb E}[L(t)]<\infty$.
\end{lemma}

\begin{remark}
The result is a special case of Proposition 6.1 of \cite{KR14}.  
We give a short proof in our simpler setting.
\end{remark}\ 

\begin{proof}
We can obtain $D$ as $D=\{x:\psi_D(x)>0\}$ where $\psi_D$ is $C^2$ on 
${\Bbb R}^d$ and $\nabla\psi_D(x)=\epsilon (x)n_D(x)$, $x\in\partial 
D$, where $\inf_{x\in\partial D}\epsilon (x)>0$.  
Observe that
\begin{eqnarray}
\psi_D(X(t))=\psi_D(X(0))+\int_0^t\nabla\psi_D(X(s))^T\sigma (X(s
))dB(s)+\int_0^t{\Bbb L}\psi_D(X(s))ds\label{itoD}\\
+\int_0^t\nabla\psi_D(X(s))\cdot\eta (X(s))dL(s).\ \nonumber\end{eqnarray}
Every term in (\ref{itoD}) has finite expectation, except 
possibly the last, but then the last must also.  Since 
$\kappa_1=\inf_{x\in\partial D}\nabla\psi_D(x)\cdot\eta (x)>0$, we have 
\[{\Bbb E}[L(t)]<\frac 1{\kappa_1}{\Bbb E}[\int_0^t\nabla\psi_D(X
(s))\cdot\eta (X(s))dL(s)]<\infty .\]

Setting $q_{\epsilon}(z)=\int_0^z\int_0^y{\bf 1}_{[0,\epsilon ]}(
u)dudy$,
\begin{eqnarray}
q_{\epsilon}(\psi_D(X(t)))&=&q_{\epsilon}(\psi_D(X(0)))+\int_0q_{
\epsilon}'(\psi_D(X(s)))\nabla\psi_D(X(s))^T\sigma (X(s))dB(s)\non\\
&&\qquad +\int_0^tq'_{\epsilon}(\psi_D(X(s))){\Bbb L}\psi_D(X(s))
ds\label{qeps}\\
&&\qquad +\int_0^tq'_{\epsilon}(\psi_D(X(s)))\nabla\psi_D(X(s))\cdot
\eta (X(s))dL(s)\non\\
&&\qquad +\int_0^t{\bf 1}_{[0,\epsilon ]}(\psi_D(X(s)))\nabla\psi_
D(X(s))^T\sigma (X(s))\sigma^T(X(s))\nabla\psi (X(s))ds.\non\end{eqnarray}
Since $\kappa_2=\inf_{x\in\bar {D}}\nabla\psi_D(x)^T\sigma (x)\sigma^
T(x)\nabla\psi (x)>0$, $\int_0^t{\bf 1}_{\partial D}(X(s))ds$ 
is bounded by $\kappa_2^{-1}$ times the last term on the right of 
(\ref{qeps}).  Since every term in (\ref{qeps}) converges 
to zero except, possibly, the last term, the last term 
also converges to zero giving the lemma.
\end{proof}

\begin{lemma}\label{bndryreg}
Let $X^x$ be as in the proof of Lemma \ref{stcsg}, and let 
$\gamma^x=\inf\{t>0:X^x(t)\in\partial D\}$.  Then 
 each $x_0\in\partial D$ is regular in 
the sense that 
\[\lim_{x\rightarrow x_0}{\Bbb E}[\gamma^x]=0\mbox{\rm \ and }\lim_{
x\rightarrow x_0}{\Bbb E}[|X^x(\gamma^x)-x_0|]=0.\]
 In particular,
\[{\Bbb E}[|X^{x_0}(\gamma^{x_0})-x_0|]={\Bbb E}[\gamma^{x_0}]=0.\]
\end{lemma}

\begin{proof}
With reference to Section 6.2 of \cite{Fri75}, a function 
$w_{x_0}\in C^2(\bar {D})$ is a 
{\em barrier\/} at $x_0$ if $w_{x_0}\geq 0$, $w_{x_0}(x)=0$ only if $
x=x_0$, and 
${\Bbb L}w_{x_0}(x)\leq -1$, and under Condition \ref{diffcnd}, a barrier 
exists for each $x_0\in\partial D$.  By It\^o's formula, 
\[w_{x_0}(X^x(t\wedge\gamma^x))-w_{x_0}(x)-\int_0^{t\wedge\gamma^
x}{\Bbb L}w_{x_0}(X^x(s))ds.\]
Taking expectations and letting $t\rightarrow\infty$, we have
\[w_{x_0}(x)\geq {\Bbb E}[w_{x_0}(X(\gamma^x))]+{\Bbb E}[\gamma^x
].\]
Since $\lim_{x\rightarrow x_0}w_{x_0}(x)=0$, the lemma follows.
\end{proof}

\begin{lemma}
The Markov process corresponding to (\ref{rsde}) has a 
stationary distribution denoted by $\pi$, and the support of 
$\pi$ is $\!\bar {D}$.  
\end{lemma}

\begin{proof}
The compactness of $\bar {D}$ and Lemma \ref{stcsg} imply 
existence of a stationary distribution by Theorem 4.9.3 
of \cite{EK86}.  To show that $\pi$ charges every open set, 
first observe that $\pi (\partial D)=0$, since the process spends 
zero real time in $\partial D$.  Let $P_x^X$ be the distribution of the 
solution $X$ with $X(0)=x$.  It then suffices to 
show that for any $x\in D$ and any ball $B\subset D$, 
$P_x^X(\tau_B<\tau_{\partial D})>0$, where $\tau_B$ denotes the first hitting 
time of $B$ and $\tau_{\partial D}$ denotes the first hitting time of the 
boundary.  By connectedness, we can find a differentiable 
path ${\cal P}$ starting at $x$  and ending at the center of the ball 
$B$ without hitting the boundary.  Let 
$\epsilon =\inf\{|x-y|:x\in {\cal P},y\in\partial D\}$.  Since $\sigma$ is nondegenerate,
by Girsanov's theorem, we can construct a
distribution $Q_x^X$ equivalent to $P_x^X$   such that under $Q_x^
X$, 
with high probability $X$ stays within $\epsilon$ of ${\cal P}$ until 
$\tau_B$.
\end{proof}

Since we are assuming that the $X_i$ are 
independent and stationary (obtained by assuming that 
the $X_i(0)$ are independent with distribution $\pi$), 
an immediate consequence of our assumptions
 is that  $V(t)$ given by (\ref{solrep}) will be absolutely continuous 
with respect to $\pi$. 

\begin{lemma}\label{abscon}
The measure $V(t)$ is absolutely continuous with respect 
to $\pi$.
\end{lemma}

\begin{proof}
The de Finetti measure for $\{(A_i(t),X_i(t))\}$ is given by the 
regular conditional distribution $\Xi (t,da,dx)$ of $(A_1(t),X_1(
t))$ 
given the tail $\sigma$-algebra ${\cal T}=\cap_n\sigma ((A_i(t),X_
i(t)):i\geq n)$. Then 
the measure $V(t)$ can be written as
\[V(t,B)=\int_{{\Bbb R}\times\bar {D}}{\bf 1}_B(x)a\Xi (t,da,dx).\]
But the $\bar {D}$-marginal of $\Xi (t,da,dx)$ is $\pi$, so there is a 
transition function ${\cal V}(t,x,da)$ satisfying
\[V(t,B)=\int_B\int_{{\Bbb R}}a{\cal V}(t,x,da)\pi (dx),\]
and hence $V(t)$ is absolutely continuous with respect to 
$\pi$.
\end{proof}

Consequently, 
we can write 
\[V(t,dx)=v(t,x)\pi (dx),\]
with 
\begin{equation}v(t,x)=\int_{{\Bbb R}}a{\cal V}(t,x,da).\label{densty}\end{equation}
As noted above, we can assume that $v$ is 
${\cal B}[0,\infty )\times {\cal B}(\bar {D})\times\vee_t{\cal F}_
t$-measurable. 

To construct a particle representation for a solution 
of (\ref{aceq}), we still need to define our weights, $A_i$.
Set $\tau_i(t)=0\vee\sup\{s\leq t:X_i(s)\in\partial D\}$, that is $
\tau_i(t)$ is the 
most recent time that $X_i$ has been on the boundary, or if 
$X_i$ has not hit the boundary by time $t$, $\tau_i(t)=0$. Of 
course, $\tau_i(t)$ is not a stopping time; however, it is 
independent of $W$,  so the stochastic integral in the 
following equation is well-defined.  We take $A_i$ to satisfy 
\begin{eqnarray}
A_i(t)&=&g(X_i(\tau_i(t))){\bf 1}_{\{\tau_i(t)>0\}}+h(X_i(0)){\bf 1}_{
\{\tau_i(t)=0\}}+\int_{\tau_i(t)}^tb(X_i(s))ds\label{gensys}\\
&&\qquad +\int_{\tau_i(t)}^tG(v(s,X_i(s)),X_i(s))A_i(s)ds+\int_{{\Bbb U}
\times (\tau_i(t),t]}\rho (X_i(s),u)W(du,ds),\nonumber\end{eqnarray}
where $v(t,x)$  is the density with respect to $\pi$ 
determined by (\ref{densty}).

As we will see in Section \ref{secteqbc}, $A_i$ does not 
appear to be a semimartingale, but the difficulties only 
occur when $X_i$ is at the boundary.  Consequently, 
we can define the integrals in the following lemma 
directly as limits of Riemann-like sums.

\begin{lemma}
For $\varphi$ in $C^2_c(D)$, 
\begin{eqnarray}
\varphi (X_i(t))A_i(t)&=&\varphi (X_i(0))A_i(0)+\int_0^t\varphi (
X_i(s))dA_i(s)\label{prespde}\\
&&\qquad +\int_0^tA_i(s)\nabla\varphi (X_i(s))^T\sigma (X_i(s))dB_
i(s)+\int_0^t{\Bbb L}\varphi (X_i(s))A_i(s)ds\non\\
&=&\varphi (X_i(0))A_i(0)+\int_0^t\varphi (X_i(s))G(v(s,X_i(s)),X_
i(s))A_i(s)ds\non\\
&&\qquad +\int_0^t\varphi (X_i(s))b(X_i(s))ds+\int_{{\Bbb U}\times 
[0,t]}\varphi (X_i(s))\rho (X_i(s),u)W(du\times ds)\non\\
&&\qquad +\int_0^tA_i(s)\nabla\varphi (X_i(s))^T\sigma (X_i(s))dB_
i(s)+\int_0^t{\Bbb L}\varphi (X_i(s))A_i(s)ds,\non\end{eqnarray}
where for partitions $\{t_i\}$ of $[0,t]$,
\begin{equation}\int_0^t\varphi (X_i(s))dA_i(s)=\lim_{\max_k|t_{k
+1}-t_k|\rightarrow 0}\sum_k\varphi (X_i(t_k))(A_i(t_{k+1})-A_i(t_
k)).\label{phiadef}\end{equation}
\end{lemma}

\begin{remark}
Since $\varphi$ vanishes in a neighborhood of the  $\partial D$, the local 
time integral in (\ref{itophi}) does not appear in this 
identity.
\end{remark}

\begin{proof}
The convergence of the limit in (\ref{phiadef}) follows by 
observing that for $\max_k|t_{k+1}-t_k|$ sufficiently small, 
 $t_k\leq\tau_i(s)\leq t_{k+1}$ implies $\varphi (X_i(t_k))=0$.  We also note that 
by this observation and the fact that the 
independence of $B_i$ and $W$ implies that the covariation of 
the two stochastic integral terms is zero,
\[\lim_{\max_k|t_{k+1}-t_k|\rightarrow 0}\sum_k(\varphi (X_i(t_{k
+1}))-\varphi (X_i(t_k)))(A_i(t_{k+1})-A_i(t_k))=0.\]
\end{proof}

To see that the weights $A_i$ determined by (\ref{gensys}) 
should give a solution of (\ref{aceq}), we want to average 
(\ref{prespde}).
The next to the last term in the right side of 
(\ref{prespde}) 
is a martingale, and these 
martingales are orthogonal for different values of $i$. 
Consequently, they will average to zero.
Assuming exchangeability of $\{(A_i,X_i)\}$, which will follow from the 
exchangeability of $\{(A_i(0),X_i)\}$ 
provided we can show existence and uniqueness for the system 
(\ref{gensys}), averaging gives (\ref{aceq}). 

We also note that the weights, at least plausibly, capture 
the desired boundary conditions.
Intuitively, for $x$ close to the boundary $\partial D$, the value of 
$v(t,x)$ is determined by the particles with locations $X_i(t)$  
close to the boundary, but if $X_i(t)$ is close to the 
boundary it should have recently hit the boundary and 
$A_i(t)$ should be close to $g(X_i(\tau_i(t)))$. 

That intuition leads 
to the following interpretation of the boundary condition.
Assume that $g$ is continuous, and let
  $\bar {g}:D\cup\partial D\mapsto\mathbb{R}$ be a continuous function 
such that $\bar {g}|_{\partial D}=g$.  
For $\epsilon >0$, define 
\[\partial_{\epsilon}D=\{x\in\bar {D}|\mathrm{dist}(x,\partial D)
<\epsilon \}.\]
Then, 
\[\int_{\partial_{\epsilon}D}|v(t,x)-\bar {g}(x)|\pi (dx)=\lim_{n
\rightarrow\infty}\frac 1n\sum_{i=1}^n{\bf 1}_{\partial_{\epsilon}
D}(X_i(t)|v(t,X_i(t))-\bar {g}(X_i(t))|.\]
Once we have existence for (\ref{gensys}), it will follow 
from Lemma \ref{cndA}\ that ${\Bbb E}[A_i(t)|W,X_i(t)]=v(t,X_i(t)
)$,
and by the intuition and the continuity of $\bar {g}$, if $X_i(t)$ is 
close to the boundary, we should 
have $v(t,X_i(t))\approx\bar {g}(X_i(t))$ and 
\begin{equation}\lim_{m\to\infty}\frac 1{\pi (\partial_{\epsilon}
D)}\int_{\partial_{\epsilon}D}|v(t,x)-\bar {g}(x)|\pi (dx)=0.\label{actualbc}\end{equation}
The intuition is made precise under regularity conditions 
on the time-reversal of the $X_i$. See Lemma \ref{bndval}
for details.  
  
Our first step will be to prove uniqueness for the 
system (\ref{gensys}) under the following condition 
which we assume throughout the paper.  

\begin{condition}\ \label{fndcndI}
The coefficients in (\ref{gensys}) satisfy
\begin{enumerate}
\item$g$ and $h$ are bounded with sup norms $\Vert g\Vert$ and $\Vert 
h\Vert$.

\item$K_1\equiv\sup_{x\in\bar {D}}|b(x)|<\infty$.

\item$K_2\equiv\sup_{x\in\bar {D}}\int\rho (x,u)^2\mu (du)<\infty$.

\item$\label{gbnd}K_3\equiv\sup_{v\in {\Bbb R},x\in\bar {D}}G(v,x
)<\infty$.

\item$L_1\equiv\sup_{v\in {\Bbb R},x\in\bar {D}}\frac {|G(v,x)|}{
1+|v|^2}<\infty$.

\item$L_2\equiv\sup_{v_1\neq v_2\in {\Bbb R},x\in\bar {D}}\frac {
|G(v_1,x)-G(v_2,x)|}{|v_1-v_2|(1+|v_1|+|v_2|)}<\infty$.  
\end{enumerate}
\end{condition}
Observe that Condition \ref{fndcndI}.\ref{gbnd}  does not 
imply that $G$ has a lower bound, but only an upper 
bound.  
For example, $G(v,x)=1-v^2$ gives the classical Allen-Cahn 
equation, whilst $G(v,x)=-v^2$ gives the $\Phi^4_d$ equation.

\begin{theorem}\label{uniqueAI}
The solution of 
(\ref{gensys}) with $v(t,x)$ the density of $V$ given by 
(\ref{solrep}) exists and is unique.
\end{theorem}

\begin{proof}
Uniqueness is proved in Section \ref{isunique} and existence 
in Section \ref{isexist}.
\end{proof}

Theorem \ref{uniqueAI}\ ensures the existence of a (signed) 
measure-valued process satisfying (\ref{aceq}) for 
$\varphi\in C_c^2(D)$.  Unfortunately, even coupled with some 
interpretation of the boundary condition, (\ref{aceq}) 
with this space of test functions does not, in general, uniquely 
determine a measure-valued process.\footnote{For example, consider $X_i$ reflecting 
Brownian motion with differing directions of reflection 
but whose stationary distribution is still normalized 
Lebesgue measure.} 
 Consequently, we 
need to enlarge the space of test functions.  We have 
two ways of doing that, first 
by taking the test functions to be $C_0^2(\bar {D})$, the space of 
twice continuously differentiable functions that vanish 
on the boundary, Theorem \ref{spdeunqI}, and second 
by taking the test functions to be 
${\cal D}({\Bbb A})$, the domain of the generator for the semigroup 
$\{{\Bbb T}(t)\}$ corresponding to the $X_i$, Theorem \ref{spdeunq2I}.

We need to identify the 
space in which the solution will live.

\begin{definition}\label{compdef}
A process $Z$ is {\em compatible\/} with a process $Y$ if
\[{\Bbb E}[f(Y)|{\cal F}_t^{Y,Z}]={\Bbb E}[f(Y)|{\cal F}_t^Y],\]
for all bounded, measurable $f$ defined on the range of $Y$ and all $
t$. 
\end{definition}

\begin{remark}\label{compind}
If $Y$ is a process with independent increments, then $Z$ is 
compatible with $Y$ if $Y(t+\cdot )-Y(t)$ is independent of ${\cal F}_
t^{Y,Z}$. 
See \cite{Kur14}, Lemma 2.4.
\end{remark}

Let ${\cal L}(\pi )$ be the space of processes $v$ compatible with $
W$ 
taking values in $L^1(\pi )$ such that for each $T>0$ and some 
$\varepsilon_T>0$, $v$ satisfies 
\[\sup_{t\leq T}{\Bbb E}\left[\int_{\bar {D}}e^{\varepsilon_T|v(t
,x)|^2}\pi (dx)\right]<\infty .\]

We prove the following theorem in Section \ref{C02}.

\begin{theorem}\label{spdeunqI} 
Consider the equation
\begin{eqnarray}
\langle\varphi (\cdot ,t),V(t)\rangle&=&\langle\varphi (\cdot ,0)
,h\rangle_{\pi}+\int_0^t\langle\varphi (\cdot ,s)G(v(s,\cdot ),\cdot 
),V(s)\rangle ds\non\\
&&\qquad +\int_0^t\int_{\bar {D}}\varphi (x,s)b(x)\pi (dx)ds\label{weakeqbca1}\\
&&\qquad +\int_{{\Bbb U}\times [0,t]}\int_{\bar {D}}\varphi (x,s)
\rho (x,u)\pi (dx)W(du\times ds)\non\\
&&\qquad +\int_0^t\langle {\Bbb L}\varphi (\cdot ,s)+\partial\varphi 
(\cdot ,s),V(s)\rangle ds\non\\
&&\qquad +\int_0^t\int_{\partial D}g(x)\eta (x)\cdot\nabla\varphi 
(x,s)\beta (dx)ds,\non\end{eqnarray}
where the  test functions $\varphi (x,t)$ are twice 
differentiable in $x$, differentiable in  $t$, and vanish on 
$\partial D\times [0,\infty )$,
$\pi$ is the stationary distribution for the particle 
location process, and $\beta$ is the measure associated with 
the local time defined in Section 
\ref{sectbndrymeas}.

Suppose ${\cal D}_0=\{\varphi\in C^2(\bar {D}):\varphi (x)=0\mbox{\rm \ and }
{\Bbb L}\varphi (x)=0,x\in\partial D\}$ is 
a core for ${\Bbb A}^0$, the generator of the semigroup $\{{\Bbb T}^
0(t)\}$ for 
the diffusion that absorbs at $\partial D$, that is, ${\Bbb A}^0$ is the 
closure of $\{(\varphi ,{\Bbb L}\varphi ):\varphi\in {\cal D}_0\}$.  Then  $
V$ defined in 
(\ref{solrep})  with $\{(X_i,A_i)\}$ given by (\ref{infsys2}) and 
(\ref{gensys}) is the unique solution of (\ref{weakeqbca1}) 
in ${\cal L}(\pi )$.  
\end{theorem}

\begin{remark}
Theorem 8.1.4 of \cite{EK86} gives conditions implying ${\cal D}_
0$ 
is a core.  Note that any solution of (\ref{weakeqbca1}) 
is a solution of (\ref{aceq}).
\end{remark}

Now we take the test functions to be ${\cal D}({\Bbb A})$, the domain 
of the generator for the semigroup $\{{\Bbb T}(t)\}$ corresponding 
to the location processes.  More precisely, let ${\cal T}$ be the 
collection of functions  $\varphi (x,t)$ for which there exists 
$t_{\varphi}>0$ such that $\varphi (t,x)=0$ for $t\ge t_{\varphi}$, $
\varphi$ is continuously 
differentiable in $t$, and $\varphi (\cdot ,t)\in {\cal D}({\Bbb A}
)$, $t\geq 0$, with ${\Bbb A}\varphi$  
bounded and continuous.

Let 
\begin{equation}\gamma_i(s)=\inf\{t>s:X_i(t)\in\partial D\},\label{gamdef1}\end{equation}
and note that 
${\bf 1}_{\{\tau_i(t)=0\}}={\bf 1}_{\{\gamma_i(0)>t\}}$.  
Let $P(dy,ds|x)$ be the conditional distribution of 
$(X_i(\gamma_i(0)),\gamma_i(0))$ given $X_i(0)=x$, and let 
$P\varphi (x)=\int\varphi (y,s)P(dy,ds|x)$. Let $X^{*}$ be the reversed 
process and $\gamma^{*}$ be the first time that $X^{*}$ hits the 
boundary.  

To simplify notation in the equation, we extend $g$ to all 
of $\bar {D}$ by setting $g(x)=h(x)$ for $x\in D$.  We do not 
require that this extension be continuous.

We prove the following theorem in Section \ref{A}.

\begin{theorem}\label{spdeunq2I} 
For the equation
\begin{eqnarray}
0&=&\int_0^{\infty}\langle (\varphi (\cdot ,s)-P\varphi (\cdot ,\cdot 
+s))G(v(s,\cdot ),\cdot ),V(s)\rangle ds\label{wkeqbc1}\\
&&\qquad +\int_0^{\infty}\int b(x)(\varphi (x,s)-P\varphi (x,\cdot 
+s))\pi (dx)ds\non\\
&&\qquad +\int_{{\Bbb U}\times [0,\infty )}\int_{\bar {D}}(\varphi 
(x,s)-P\varphi (x,\cdot +s))\rho (x,u)\pi (dx)W(du\times ds)\non\\
&&\qquad +\int_0^{\infty}\langle {\Bbb A}\varphi (\cdot ,s)+\partial
\varphi (\cdot ,s),V(s)\rangle ds\non\\
&&\qquad -\int_0^{\infty}\int_{\bar {D}}{\Bbb E}[g(X^{*}(\gamma^{
*}\wedge s))|X^{*}(0)=x]({\Bbb A}\varphi (x,s)+\partial\varphi (x
,s))\pi (dx)ds.\non\end{eqnarray}
for all $\varphi\in {\cal T}$, 
 $V$ defined in 
(\ref{solrep}) with $\{(X_i,A_i)\}$ given by (\ref{infsys2}) and 
(\ref{gensys}) is the unique solution of (\ref{wkeqbc1}) in 
${\cal L}(\pi )$.  
\end{theorem}

\begin{remark}
The form of (\ref{wkeqbc1}) may not be very intuitive; 
however, this equation and (\ref{weakeqbca1}) determine 
the same unique solution in ${\cal L}(\pi )$.  

To see that the 
restriction of (\ref{wkeqbc1}) to $\varphi\in C_c^2(D)$ gives 
(\ref{aceq}), note first that for $\varphi\in C_c^2(D)$, $P\varphi 
=0$.  To 
see that the last term in (\ref{wkeqbc1}) is zero for 
$\varphi\in C_c^2(D)$, note that $u(x,s)={\Bbb E}[g(X^{*}(\gamma^{
*}\wedge s))|X^{*}(0)=x]$ 
should be a solution of the Dirichlet problem 
$({\Bbb A}_0^{*}-\partial )u(x,s)=0$ on $\bar {D}\times [0,\infty 
)$ with boundary conditions 
$u(x,0)=h(x)$, $x\in D$, and $u(x,s)=g(x)$, $s>0$, $x\in\partial 
D$.  
\end{remark}

\setcounter{equation}{0}

\section{Existence and uniqueness of the weighted 
particle system}\label{sectaddnoise} 
In this section we prove a number of estimates for particle systems of the type we consider in this paper, leading up to a proof of Theorem \ref{uniqueAI}. The system of stochastic differential equations
 (\ref{gensys}) must be considered in 
conjunction with the existence of an empirical 
distribution 
\[V(t)=\lim_{n\rightarrow\infty}\frac 1n\sum_{i=1}^nA_i(t)\delta_{
X_i(t)}\]
required to have a density $v(t,\cdot )$ with respect to $\pi$.  It 
is by no means clear that a solution satisfying all these 
constraints exists.  

First we explore the properties that a solution must 
have by replacing $v$ by an arbitrary, measurable  
$L^1(\pi )$-valued stochastic process  $U$ that is independent of 
$\{X_i\}$ and compatible with $W$, (see Definition 
\ref{compdef}).  In the current setting, compatibility 
means that for each $t>0$, 
$\sigma (W(C\times (t,t+s]):C\in {\cal B}({\Bbb U}),\mu (C)<\infty 
,s>0)$ is independent of 
\[{\cal F}^{U,W}_t=\sigma (U(s),W(C\times [0,s]):0\leq s\leq t,C\in 
{\cal B}({\Bbb U}),\mu (C)<\infty ).\]

Define $A_i^U$ 
to be the solution of
\begin{eqnarray}
A^U_i(t)&=&g(X_i(\tau_i(t))){\bf 1}_{\{\tau_i(t)>0\}}+h(X_i(0)){\bf 1}_{
\{\tau_i(t)=0\}}\label{avdef}\\
&&\qquad +\int_{\tau_i(t)}^tG(U(s,X_i(s)),X_i(s))A^U_i(s)ds+\int_{
\tau_i(t)}^tb(X_i(s))ds\nonumber\\
&&\qquad +\int_{{\Bbb U}\times (\tau_i(t),t]}\rho (X_i(s),u)W(du\times 
ds).\nonumber\end{eqnarray}
Existence and uniqueness of the solution of (\ref{avdef}) 
holds under modest assumptions on the coefficients, in 
particular, under Condition \ref{fndcndI}.

\begin{lemma}\label{solest}
Let
\[H_i(t)=\int_{{\Bbb U}\times [0,t]}\rho (X_i(s),u)W(du,ds).\]
Then $H_i$ is a martingale with respect to the filtration 
$\{{\cal H}^i_t\}\equiv \{{\cal F}_t^W\vee\sigma (X_i)\}$ and there exists a standard Brownian motion 
$Z_i$ such that $Z_i$ is independent of $X_i$ and
\[H_i(t)=Z_i\left(\int_0^t\rho^2(X_i(s),u)\mu (du)\right).\]
For all $A_i^U$ defined as in (\ref{avdef}) and $K_1$, $K_2$, and $
K_3$ 
defined in Condition \ref{fndcndI},
\begin{eqnarray}
|A^U_i(t)|&\leq&(\Vert g\Vert\vee\Vert h\Vert +K_1(t-\tau_i(t))+\sup_{
\tau_i(t)\leq r\leq t}|H_i(t)-H_i(r)|)e^{K_3(t-\tau_i(t))}
\label{abnd}\\
&\leq&(\Vert g\Vert\vee\Vert h\Vert +K_1t+\sup_{0\leq s\leq t}|H_
i(t)-H_i(s)|)e^{K_3t}\equiv\bar {A}_i(t)\nonumber\\
&\leq&(\Vert g\Vert\vee\Vert h\Vert +K_1t+2\sup_{0\leq s\leq t}|Z_
i(sK_2)|)e^{K_3t}\equiv\Gamma_i(t).\nonumber\end{eqnarray}
For each $T>0$, there exists $\varepsilon_T$ such that 
\begin{equation}{\Bbb E}[e^{\varepsilon_T\sup_{t\leq T}|A_i^U(t)|^
2}]\leq {\Bbb E}[e^{\varepsilon_T\Gamma_i(T)^2}]<\infty .\label{expbnd}\end{equation}
\end{lemma}

\begin{proof}
Let $A^{+}_i(t)=A_i^U(t)\vee 0$ and $A^{-}_i(t)=(-A_i^U(t))\vee 0$. Define 
\[\gamma_i^{+}(t)=\tau_i(t)\vee\sup\{s<t:A_i^U(s)<0\}.\]
If $\gamma_i^{+}(t)=\tau_i(t)<t$, then $A_i^U(s)\geq 0$ for all $
s$ in 
(\ref{avdef});
if $0<\gamma_i^{+}(t)<t$, $A_i^U(\gamma_i^{+}(t))=0$ and  
\begin{eqnarray*}
A^U_i(t)&=&\int_{\tau_i(t)}^tG(U(s,X_i(s)),X_i(s))A^U_i(s)ds+\int_{
\tau_i(t)}^tb(X_i(s))ds\\
&&\qquad +\int_{{\Bbb U}\times (\tau_i(t),t]}\rho (X_i(s),u)W(du\times 
ds);\end{eqnarray*}
if $\gamma_i^{+}(t)=t$, $A_i^{+}(t)\leq\Vert g\Vert\vee\Vert h\Vert$.  In any of these cases
\[A^{+}_i(t)\leq\Vert g\Vert\vee\Vert h\Vert +\int_{\gamma_i^{+}(
t)}^tK_3A^{+}_i(s)ds+K_1(t-\gamma_i^{+}(t))+\sup_{\gamma_i^{+}(t)
\leq s\leq t}|H_i(t)-H_i(s)|,\]
so by Gronwall,
\[A^{+}_i(t)\leq (\Vert g\Vert\vee\Vert h\Vert +K_1(t-\gamma_i^{+}
(t))+\sup_{\gamma_i^{+}(t)\leq s\leq t}|H_i(t)-H_i(s)|)e^{K_3(t-\gamma^{
+}_i(t))}.\]
Letting $\gamma^{-}_i(t)=\tau_i(t)\vee\sup\{s<t:A_i^U(s)>0\}$, similar 
observations give
\[A^{-}_i(t)\leq\Vert g\Vert\vee\Vert h\Vert +\int_{\gamma^{-}_i(
t)}^tK_3A^{-}_i(s)ds+K_1(t-\gamma^{-}_i(t))+\sup_{\gamma_i^{-}(t)
\leq s\leq t}|H_i(t)-H_i(s)|,\]
so we have a similar bound on $A^{-}_i$.  Together the bounds 
give the first two inequalities in (\ref{abnd}).

$H_i$ is a continuous 
martingale with quadratic variation 
$\int_0^t\int\rho (X_i(s),u)^2\mu (du)ds$.  Define
\[\gamma (u)=\inf\{t:\int_0^t\int\rho (X_i(s),u)^2\mu (du)ds\geq 
u\},\]
and $Z_i(u)=H_i(\gamma (u))$.  Then $Z_i$ is a continuous martingale 
with respect to the filtration $\{{\cal H}^i_{\gamma (u)}\}$ and $
[Z_i]_u=u$, so $Z_i$ 
is a standard 
Brownian motion. Since $\sigma (X_i)\subset {\cal H}^i_0$, $Z_i$ is independent 
of $X_i$.  

The first inequality in (\ref{expbnd}) follows by the 
monotonicity of $\Gamma_i$ and the finiteness by 
standard estimates on the distribution of the 
supremum of Brownian motion.
\end{proof}

The pairs $\{(A_i^U,X_i)\}$ will be exchangeable, so with reference to 
Lemma \ref{abscon}, we can define 
$\Phi U(t,x)$ to be the density with respect to $\pi$ of the signed measure 
determined by 
\begin{equation}\langle\varphi ,\Phi U(t)\rangle =\lim_{N\rightarrow
\infty}\frac 1N\sum_{i=1}^NA_i^U(t)\varphi (X_i(t)).\label{Phidef}\end{equation}

\begin{lemma}\label{cndA}
Suppose that $U$ is compatible with $W$ and that $(U,W)$ is 
independent of $\{X_i\}$.  Then $\Phi U$ is compatible with $W$
and
for each $i$,
\begin{equation}{\Bbb E}[A_i^U(t)|U,W,X_i(t)]=\Phi U(t,X_i(t)).\label{cndrep}\end{equation}
If, moreover, $U$ is $\{{\cal F}_t^W\}$-adapted, then $\Phi U$ is $
\{{\cal F}_t^W\}$-adapted.
\end{lemma}

\begin{proof}
It follows from (\ref{Phidef}) that $\Phi U(t)$ is measurable 
with respect to the shift invariant sigma algebra of the 
stationary sequence $(X_i,U,W)$.  Since the $\{X_i\}$ sequence is 
i.i.d.  and independent of $(U,W)$, this sigma algebra is 
contained in the completion of the sigma algebra 
generated by $(U,W)$.  It then follows from the ergodic 
theorem that 
\begin{eqnarray*}
\langle\varphi ,\Phi U(t)\rangle ={\Bbb E}[A_1^U(t)\varphi (X_1(t
))|U,W].\end{eqnarray*}
Since $A_i^U$ is $\{{\cal F}_t^{X_i,U,W}\}$-adapted, $X_i$ is $\{
{\cal F}_t^{X_i}\}$-adapted, and 
$X_i$ and $(U,W)$ are independent, we may replace 
conditioning by $(U,W)$ by conditioning by ${\cal F}_t^{U,W}$ in the 
previous expression, which shows the compatibility for 
$\Phi U$ in general and $\{{\cal F}_t^W\}$-adaptedness if $U$ is $
\{{\cal F}_t^W\}$-adapted.  
By exchangeability, 
\begin{eqnarray*}
{\Bbb E}[A_i^U(t)\varphi (X_i(t))F(U,W)]&=&{\Bbb E}[\int\varphi (
x)\Phi U(t,x)\pi (dx)F(U,W)]\\
&=&{\Bbb E}[\varphi (X_i(t))\Phi U(t,X_i(t))F(U,W)],\end{eqnarray*}
where the second equality follows by the independence of 
$X_i(t)$ and $(U,W)$.  The lemma then follows by the definition 
of conditional expectation.
\end{proof}

We will restrict attention in the following results to 
the case where $U$ is $\{{\cal F}_t^W\}$-adapted to simplify the 
notation slightly.  Analogues of these results hold for 
compatible $U$ as well.  The next result shows existence 
of a version of the density $\Phi U(t,x)$ with the property 
that $t\mapsto\Phi U(t,X_i(t))$ is well-behaved pathwise.  
\begin{lemma}\label{filtcnd} 
Suppose that $U$ is $\{{\cal F}_t^W\}$-adapted, and let 
${\cal G}_t^{X_i}=\sigma (X_i(r):r\geq t)$.  Then there exists a version of $
\Phi U(t,x)$ 
such that 
\begin{equation}\Phi U(t,X_i(t))={\Bbb E}[A_i^U(t)|W,X_i(t)]={\Bbb E}
[A_i^U(t)|\sigma (W)\vee {\cal G}_t^{X_i}],\label{revmp}\end{equation}
where we interpret the right side as the optional 
projection, and for this version
\begin{equation}{\Bbb E}[\sup_{0\leq t\leq T}|\Phi U(t,X_i(t))|^2
]\leq 4{\Bbb E}[\sup_{0\leq t\leq T}|A_i^U(t)|^2].\label{phidoob}\end{equation}

Moreover the identity (\ref{cndrep}) holds with $t$ replaced by any 
nonnegative $\sigma (W)$-measurable random variable $\tau$.
\end{lemma}

\begin{proof}
The first equality in (\ref{revmp}) is just (\ref{cndrep}), 
and the second follows from the fact that $X_i$ is Markov.
By Lemma \ref{opmeas}, there exists a Borel measurable 
function $g$ on $[0,\infty )\times {\Bbb R}^d\times {\Bbb W}$ such that
\[g(t,X_i(t),W)={\Bbb E}[A_i^U(t)|\sigma (W)\vee {\cal G}_t^{X_i}
].\]
It follows that $g(t,x,W)$ is a version of 
$\Phi U(t,x)$. (Note that an arbitrary version of $\Phi U(t,x)$ may 
not have the measurability properties of $g(t,x,W)$.)

Corollary \ref{opmon}, the properties of reverse 
martingales, and Doob's inequality give (\ref{phidoob}),  
and the last statement follows by the definition of the 
optional projection.
\end{proof}

With (\ref{cndrep}) in mind, given an exchangeable family 
$\{A_i\}$ such that $A_i$ is adapted to $\{{\cal F}^{X_i}_t\vee {\cal F}^W\}$, 
define $\Phi A_i\equiv A^U_i$ taking $U$ to be given by
\[U(t,X_i(t))={\Bbb E}[A_i(t)|W,X_i(t)]={\Bbb E}[A_i(t)|\sigma (W
)\vee {\cal G}^{X_i}_t]\]
in (\ref{avdef}).

\begin{lemma}\label{condbd} 
Suppose that $U$ is $\{{\cal F}_t^W\}$-adapted. Then for all $T>0$ and $
p\geq 1$, there exists a constant $C_{p,T}$ so 
that for all $t\leq T$ 
\begin{eqnarray*}
{\Bbb E}[|\Phi U(t,X_i(t))|^p|\mathcal{G}_T^{X_i}]\leq C_{p,T}\end{eqnarray*}
and 
\begin{eqnarray*}
{\Bbb E}[|A_i^U(t)|^p|\mathcal{G}_T^{X_i}]\leq C_{p,T}.\end{eqnarray*}
\end{lemma}
\begin{proof}
Recall that we have the bound
\begin{eqnarray*}
|A_i^U(t)|\leq\left(\|g\|\vee\|h\|+K_1t+2\sup_{0\leq r\leq t}|\int_{\mathbb{
U}\times (0,r]}\rho (X_i(s),u)W(du\times ds)|\right)e^{K_3t}\end{eqnarray*}
and that
\begin{eqnarray*}
\Phi U(t,X_i(t))={\Bbb E}[A_i^U(t)|W,\mathcal{G}^{X_i}_t].\end{eqnarray*}
Notice that Jensen's inequality gives
\begin{eqnarray*}
{\Bbb E}[|\Phi U(t,X_i(t))|^p|\mathcal{G}_T^{X_i}]&=&{\Bbb E}\left
[\Big|{\Bbb E}[A_i^U(t)|W,\mathcal{G}_t^{X_i}]\Big|^p\bigg|\mathcal{
G}_T^{X_i}\right]\\
&\leq&{\Bbb E}\left[{\Bbb E}\left[|A_i^U(t)|^p\Big|W,\mathcal{G}_
t^{X_i}\right]\Big|\mathcal{G}_T^{X_i}\right]\end{eqnarray*}
and that $t\leq T$ implies 
$\mathcal{G}_T^{X_i} \subset \mathcal{G}_t^{X_i} \vee \sigma(W)$. It follows that
\begin{eqnarray*}
{\Bbb E}\left[{\Bbb E}\left[|A_i^U(t)|^p\Big|W,\mathcal{G}_t^{X_i}\right
]\Big|\mathcal{G}_T^{X_i}\right]&=&{\Bbb E}\left[|A_i^U(t)|^p\Big
|\mathcal{G}_T^{X_i}\right]\leq {\Bbb E}\left[\left.\Gamma_i^p(t)\right
|{\cal G}_T^{X_i}\right].\end{eqnarray*}
Fix $S\in\mathcal{G}_T^{X_i}$ with $P(S)>0$, so that $W$ is an 
$\{\mathcal{F}_t^W\vee\sigma (X_i)\}$-martingale measure under $P
(\cdot |S)$. By the Burkholder-Davis-Gundy inequality, we find
\begin{eqnarray*}
&&{\Bbb E}\left[\Big|\sup_{0\leq r\leq t}|\int_{\mathbb{U}\times 
(0,r]}\rho (X_i(s),u)W(du\times ds)\Big|^p\bigg|S\right]\\
&&\qquad\leq C_p{\Bbb E}\left
[\left(\int_0^t\int_{\mathbb{U}}\rho (X_i(s),u)^2\mu (du)ds\right
)^{\frac p2}\Big|S\right]\\
&&\qquad\leq C_p(K_2t)^{\frac p2}.\end{eqnarray*}
$S$ was arbitrary, so the result follows.
\end{proof}

\begin{lemma}\ \label{opest}
Suppose that $U$ is $\{{\cal F}_t^W\}$-adapted. Then for each $T\geq 
0$, there exists
 $\varepsilon_T>0$ such that 
\begin{equation}{\Bbb E}\left[e^{\varepsilon_T\sup_{t\leq T}|\Phi 
U(t,X_i(t))|^2}\right]<\infty .\label{expmom}\end{equation}
\end{lemma}

\begin{proof}
As in Lemma \ref{solest}, for each $T>0$, 
there exists $\varepsilon_T>0$  such that
${\Bbb E}[e^{\varepsilon_T\Gamma_i(T)^2}]<\infty$. Recalling that 
\[|\Phi U(t,X_i(t))|=|{\Bbb E}[A_i^U(t)|\sigma (W)\vee {\cal G}_t^{
X_i}]|\leq {\Bbb E}[\Gamma_i(T)|\sigma (W)\vee {\cal G}_t^{X_i}],\]
Jensen's and Doob's inequalities give
\begin{eqnarray*}
{\Bbb E}[e^{\varepsilon_T\sup_{t\leq T}|\Phi U(t,X_i(t)|^2}]&\leq&
{\Bbb E}[\sup_{t\leq T}{\Bbb E}[e^{\varepsilon_T\Gamma_i(T)^2}|\sigma 
(W)\vee {\cal G}_t^{X_i}]]\\
&\leq&4{\Bbb E}[e^{\varepsilon_T\Gamma_i(T)^2}].\end{eqnarray*}
\end{proof}

\subsection{Proof of uniqueness in Theorem \ref{uniqueAI}}\label{isunique}
Truncations based on the moment estimates given above 
allow us to apply Gronwall's inequality to prove 
uniqueness for the system (\ref{gensys}).  

Let $U_1$ and $U_2$ 
satisfy
 $|U_k(t,X_i(t))|\leq {\Bbb E}[\Gamma_i(t)|W,X_i(t)]$, $k=1,2$.  
Then there 
exists a constant $L_3>0$ such that
\begin{eqnarray*}
&&|A_i^{U_1}(t)-A_i^{U_2}(t)|\\
&&\qquad\qquad\leq\int_{\tau_i(t)}^t|G(U_1(s,X_i(s)),X_i(s))A^{U_
1}_i(s)-G(U_2(s,X_i(s)),X_i(s))A_i^{U_2}(s)|ds\\
&&\qquad\qquad\leq\int_{\tau_i(t)}^tL_1(1+{\Bbb E}[\Gamma_i(s)|W,
X_i(s)]^2)|A_i^{U_1}(s)-A_i^{U_2}(s)|ds\\
&&\qquad\qquad\qquad +\int_{\tau_i(t)}^tL_2(1+2{\Bbb E}[\Gamma_i(
s)|W,X_i(s)]\Gamma_i(s))|U_1(s,X_i(s))-U_2(s,X_i(s))|ds\\
&&\qquad\qquad\leq\int_0^tL_1(1+C^2)|A_i^{U_1}(s)-A_i^{U_2}(s)|ds\\
&&\qquad\qquad\qquad +\int_0^tL_2(1+2C^2)|U_1(s,X_i(s))-U_2(s,X_i
(s))|ds\\
&&\qquad\qquad\qquad +\int_0^t{\bf 1}_{\{\Gamma_i(s)>C\}\cup \{{\Bbb E}
[\Gamma_i(s)|W,X_i(s)]>C\}}\Gamma_i(s)L_3(1+{\Bbb E}[\Gamma_i(s)|
W,X_i(s)]^2)ds.\end{eqnarray*}

Suppose $U_k=\Phi U_k$, $k=1,2$, that is, we have two solutions. Then 
conditioning both sides of the above inequality on $W$ and observing
\[{\Bbb E}[|U_1(s,X_i(s))-U_2(s,X_i(s))||W]\leq {\Bbb E}[|A_i^{U_
1}(s)-A_i^{U_2}(s)||W],\]
we have 
\begin{eqnarray*}
{\Bbb E}[|A_i^{U_1}(t)-A_i^{U_2}(t)||W]\leq e^{(L_1+2L_2)(1+C^2)t}
\int_0^t{\Bbb E}[{\bf 1}_{\{\Gamma_i(s)>C\}\cup \{{\Bbb E}[\Gamma_
i(s)|W,X_i(s)]>C\}}\\
\times\Gamma_i(s)L_3(1+{\Bbb E}[\Gamma_i(s)|W,X_i(s)]^2)|W]ds.\end{eqnarray*}
Taking expectations of both sides and applying H\"older's 
inequality, 
\begin{eqnarray*}
&&{\Bbb E}[|A_i^{U_1}(t)-A_i^{U_2}(t)|]\\
&&\quad\leq e^{(L_1+2L_2)(1+C^2)t}L_3\int_0^t(P\{\Gamma_i(s)>C\}^{
1/3}+P\{{\Bbb E}[\Gamma_i(s)|W,X_i(s)]>C\}^{1/3})\\
&&\qquad\qquad\qquad\qquad\qquad\qquad\times {\Bbb E}[\Gamma_i(s)^
3]^{1/3}{\Bbb E}[(1+{\Bbb E}[\Gamma_i(s)|W,X_i(s)]^2)^3]^{1/3}ds\\
&&\quad\leq e^{(L_1+2L_2)(1+C^2)t}e^{-\varepsilon_TC^2/3}2L_3\int_
0^t{\Bbb E}[e^{\frac {\varepsilon_T}3\Gamma_i(s)^2}{\Bbb E}[\Gamma_
i(s)^3]^{1/3}\\
&&\qquad\qquad\qquad\qquad\qquad\qquad\qquad\qquad\times {\Bbb E}
[(1+{\Bbb E}[\Gamma_i(s)|W,X_i(s)]^2)^3]^{1/3}ds,\end{eqnarray*}
so for $t<\frac {\varepsilon_T}{3(L_1+2L_2)}$, the right side goes to zero as 
$C\rightarrow\infty$ implying $U_1=U_2$  on $[0,t]$.  The same argument 
and induction
extends uniqueness to any interval.

\subsection{Proof of existence in Theorem \ref{uniqueAI}}\label{isexist}
The estimates of the previous section can be applied to 
give convergence of an iterative sequence proving 
existence.  For an exchangeable family $\{A_i^{(1)}\}$ with $A_i^{
(1)}$ 
$\{{\cal F}_t^{X_i,W}\}$-adapted and satisfying 
$|A_i^{(1)}(t)|\leq\Gamma_i(t)$,
recursively define $A_i^{(n+1)}=\Phi A_i^{(n)}$.  By the estimates of 
the previous section, for $n,m\geq 1$
\begin{eqnarray*}
&&|A_i^{(n+1)}(t)-A_i^{(n+m+1)}(t)|\\
&&\quad\qquad\leq\int_{\tau_i(t)}^t|G({\Bbb E}[A_i^{(n)}(s)|W,X_
i(s)],X_i(s))A^{(n+1)}_i(s)\\
&&\quad\qquad\qquad\qquad\qquad\qquad\qquad -G({\Bbb E}[A_i^{(n+
m)}(s)|W,X_i(s)],X_i(s))A_i^{(n+m+1)}(s)|ds\\
&&\quad\qquad\leq\int_{\tau_i(t)}^tL_1(1+{\Bbb E}[\Gamma_i(s)|W,
X_i(s)]^2)|A_i^{(n+1)}(s)-A_i^{(n+m+1)}(s)|ds\\
&&\quad\qquad\qquad +\int_{\tau_i(t)}^tL_2(1+2{\Bbb E}[\Gamma_i(
s)|W,X_i(s)]\Gamma_i(s))|{\Bbb E}[A_i^{(n)}(s)-A^{(n+m)}_i(s)|W,X_
i(s)]|ds\\
&&\quad\qquad\leq\int_0^tL_1(1+C^2)|A_i^{(n+1)}(s)-A_i^{(n+m+1)}
(s)|ds\\
&&\quad\qquad\qquad +\int_0^tL_2(1+2C^2)|{\Bbb E}[A_i^{(n)}(s)-A^{
(n+m)}_i(s)|W,X_i(s)]|ds\\
&&\quad\qquad\qquad +\int_0^t{\bf 1}_{\{\Gamma_i(s)>C\}\cup \{{\Bbb E}
[\Gamma_i(s)|W,X_i(s)]>C\}}\Gamma_i(s)L_3(1+{\Bbb E}[\Gamma_i(s)|
W,X_i(s)]^2)ds.\end{eqnarray*}
Setting
\[H_C(t)={\Bbb E}[{\bf 1}_{\{\Gamma_i(s)>C\}\cup \{{\Bbb E}[\Gamma_
i(s)|W,X_i(s)]>C\}}\Gamma_i(s)L_3(1+{\Bbb E}[\Gamma_i(s)|W,X_i(s)
]^2)|W],\]
we have 
\begin{eqnarray*}
&&{\Bbb E}[\sup_{r\leq t}|A_i^{(n+1)}(r)-A_i^{(n+m+1)}(r)||W]\\
&&\quad\qquad\leq\int_0^tL_1(1+C^2){\Bbb E}[\sup_{r\leq s}|A_i^{
(n+1)}(r)-A_i^{(n+m+1)}(r)||W]ds\\
&&\quad\qquad\qquad\qquad +\int_0^tL_2(1+2C^2){\Bbb E}[\sup_{r\leq 
s}|A_i^{(n)}(r)-A_i^{(n+m)}(r)||W]ds\\
&&\quad\qquad\qquad\qquad +\int_0^tH_C(s)ds.\end{eqnarray*}
It follows that
\[\limsup_{n\rightarrow\infty}\sup_{m\geq 1}{\Bbb E}[\sup_{r\leq 
t}|A_i^{(n)}(r)-A_i^{(n+m)}(r)||W]\leq e^{t(L_1+2L_2)(1+C^2)}\int_
0^tH_C(s)ds,\]
and as in the previous section, for $t<\frac {\varepsilon_T}{3(L_
1+2L_2)}$, the 
expectation of the right side goes to zero as $C$ goes to 
infinity.  Consequently, the sequence $\{A_i^{(n)}\}$ is Cauchy on 
$[0,T]$ for any fixed $T <\frac {\varepsilon_T}{3(L_
1+2L_2)}$. Then there exists $A_i$ such that 
$\lim_{n\rightarrow\infty}{\Bbb E}[\sup_{r\leq t}|A_i^{(n)}(r)-A_
i(r)|]=0$ giving existence of a 
solution on 
the interval $[0,T]$. The same argument gives existence on $[T/2,3T/2]$ and 
uniqueness shows that these solutions coincide on $[T/2,T]$. Using induction, we deduce the global 
existence of the solution of 
(\ref{gensys}).

\section{Boundary behavior} \label{sectbndrymeas} In this 
section, we present two senses in which the particle 
representation satisfies the boundary condition 
(\ref{dirchbc}).  These results depend on the boundary 
regularity of the stationary diffusions $\{X_i\}$ run forward 
or backward in time.  We begin with the result coming 
from regularity of the forward process proved in Lemma 
\ref{bndryreg}, which leads to the weak formulation of 
the stochastic PDE including the boundary condition in 
Theorem \ref{spdeunqI}. 

Let $X$ satisfy (\ref{rsde}), and assume $X(0)$ has 
distribution $\pi$ so that $X$ is stationary. 
By stationarity, for $t\in\mathbb{R}_{+}$, the process $X(t+\cdot 
)$ 
has the same distribution as $X(\cdot )$.  Therefore 
\begin{eqnarray*}
{\Bbb E}\left[\int_s^t\varphi (X(r))dL(r)\right]&=&{\Bbb E}\left
[\int_0^{t-s}\varphi (X(r))dL(r)\right].\end{eqnarray*}
For bounded and continuous $\varphi$, define
\begin{eqnarray*}
Q(t,\varphi )={\Bbb E}\left[\int_0^t\varphi (X(s))dL(s)\right
]\end{eqnarray*}
which, by the above discussion, satisfies 
$Q(t+s,\varphi )=Q(t,\varphi )+Q(s,\varphi )$ and $|Q(t,\varphi )
|\leq tC\|\varphi\|_{\infty}$ for some 
constant $C$.  Therefore, 
since $Q$ is additive in its first coordinate, 
there exists a constant $C_{\varphi}$ so that $Q(t,\varphi )=tC_{
\varphi}$.  Since $Q$ 
is also linear in its second coordinate, it then follows 
from the Riesz representation theorem that there exists 
a measure $\beta$ on $\partial D$ which satisfies 
\begin{eqnarray*}
\varphi\mapsto\frac 1t{\Bbb E}\left[\int_0^t\varphi (X(s))dL(
s)\right]&=&\int_{\partial D}\varphi (x)\beta (dx).\end{eqnarray*}

By considering test functions of product form which are 
step functions in time, we can see that for sufficiently 
regular space-time functions $\varphi$, we have 
\begin{equation}\int_0^t\int_{\partial D}\varphi (x,s)\beta (dx)d
s={\Bbb E}\left[\int_0^t\varphi (X(s),s)dL(s)\right].\label{bmeasid}\end{equation}
Denote partial derivatives with respect to the time 
variable by $\partial$.  Applying Ito's lemma to $\varphi (X(t),t
)$ for 
sufficiently smooth $\varphi$ and taking expectations, we also 
have the following relation between $\pi$ and $\beta$:  
\begin{eqnarray*}
\int_0^t\int_D\left(\partial +{\Bbb L}\right)\varphi (x,s)\pi (dx
)ds&=&\int_0^t\int_{\partial D}\nabla\varphi (x,s)\cdot\eta (x)\beta 
(dx)ds.\end{eqnarray*}

Before the next result, we recall some definitions from 
analysis.  Given a set $A\subset\mathbb{R}$ and $a\in A$, we say 
that $a$ is an {\em isolated point\/} if there exists $\epsilon >
0$ such 
that $A\cap (a-\epsilon ,a+\epsilon )=\{a\}$.  We say that $a\in 
A$ is 
{\em left-isolated\/} if there exists $\epsilon >0$ such that 
$A\cap (a-\epsilon ,a)=\emptyset$.  

\begin{lemma}\label{bdytimes} 
Almost 
surely, the set $\{t\geq 0:X(t)\in\partial D\}$ is a closed set with no 
isolated points and the collection of left-isolated points 
of this set is countable.  
\end{lemma}

The proof of Lemma \ref{bdytimes} is the same as the 
proof of the analogous property for the zero set of one 
dimensional Brownian motion and uses the results in 
Lemma \ref{bndryreg}.  See, for 
example, the proof of Theorem 2.28 in \cite{MP10}.  

\begin{lemma}\label{ltimeprop} 
Let $\tau (t)=0\vee\sup\{s\leq t:X(s)\in\partial D\}$.  Then, almost surely, 
\[\int_{\left\{t:\tau (t)\neq t\right\}}dL(s)=\int_{\left\{t:\tau 
(t-)\neq t\right\}}dL(s)=0.\]
\end{lemma}

\begin{proof}
Local time is a continuous measure supported on the set 
$\{t\geq 0:X(t)\in\partial D\}$ and therefore assigns measure zero to 
the (countable) set of left isolated points of 
$\{t\geq 0:X(t)\in\partial D\}$.  If $t_0\in \{t\geq 0:X(t)\in
\partial D\}$ is not 
left-isolated, then $t_0=\tau(t_0)=\tau(t_0-)$.  
\end{proof}

By (\ref{bmeasid}) and Lemma \ref{ltimeprop}, we have 
the following theorem.

\begin{theorem}\label{bndryid}
Almost surely, for $dL_i$ almost every $t$, 
$A_i(t)=A_i(t-)=g(X_i(t))$ and therefore 
\begin{eqnarray*}
&&\lim_{n\rightarrow\infty}\frac 1n\sum_{i=1}^n\int_0^tA_i(s-)\eta 
(X_i(s))\cdot\nabla\varphi (X_i(s),s)dL_i(s) \\
&=&{\Bbb E}\left[\int_
0^tA_i(s-)\eta (X_i(s))\cdot\nabla\varphi (X_i(s),s)dL_i(s)|U,W\right]\\
&=&\int_0^t\int_{\partial D}g(x)\eta (x)\cdot\nabla\varphi (x,s)\beta 
(dx)ds.\end{eqnarray*}
\end{theorem}

In the next section, we show how Theorem \ref{bndryid} 
leads to a weak formulation of the stochastic partial 
differential equation with a broader class of test 
functions than those considered above.  We now turn to 
the form of the boundary condition mentioned in the 
introduction at (\ref{actualbc}), which depends on 
regularity of the time-reversed process.  

For each $t$ and  $s\leq t$, define the 
time reversal of $X_i$ by $X_{i,t}^{*}(s)=X_i(t-s)$.  For notational 
convenience, when it is clear from context what the 
value of $t$ is, we will suppress the subscript and take 
the convention that $X_{i,t}^{*}(s)\equiv X_i^{*}(s)$.  Since $X_
i$ is 
stationary, the time reversal $X_i^{*}$ is a Markov process 
whose generator ${\Bbb A}^{*}$ satisfies 
\[\int_{\overline{D}}g{\Bbb A}fd\pi =\int_{\overline{D}}f{\Bbb A}^{*}gd\pi ,\quad f\in {\cal D}
({\Bbb A}).\]

\begin{remark}\ 
\label{BMExample} If $D$ is sufficiently smooth and ${\Bbb A}=\Delta$ 
with normally reflecting boundary conditions, then $\pi$ is 
proportional to Lebesgue measure, $\beta$ is proportional to 
the surface measure, and ${\Bbb A}^{*}={\Bbb A}$.  
\end{remark}

Define the hitting time of the boundary for the reversed 
process by $\sigma_i=\inf\{s:X_i^{*}(s)\in\partial D\}$, so if the reversal is 
from time $t$, $\sigma_i=t-\tau_i(t)$.  Showing that (\ref{actualbc}) 
is satisfied will depend on the following condition.  

\begin{condition}\label{bwdreg} 
The boundary $\partial D$ is regular for $X_i^{*}$ in the sense that for 
each $\delta >0$ and $x\in\partial D$, 
\begin{equation}\lim_{y\in D\rightarrow x}P\left(\sigma_i>\delta 
|X_i^{*}(0)=y\right)=0,\label{bdyhit}\end{equation}
and
\begin{equation}\lim_{y\in D\to x}{\Bbb E}\left[|X_i^{*}(\sigma_i
)-x|\wedge 1|X_i^{*}(0)=y\right]=0.\label{bdyreg}\end{equation}
\end{condition}
\begin{remark}
Condition \ref{bwdreg} is difficult to verify in general.  
A natural sufficient condition for Condition \ref{bwdreg} 
to hold is that the time-reversal is again a regular 
diffusion; see \cite{cat88} for results in this direction.  
The motivating examples discussed in the introduction all 
have $\mathbb{L}=\Delta$ and are naturally considered with 
respect to $\pi$ proportional to the Lebesgue measure, $\beta$ 
proportional to the surface measure of $\partial D$, and reflection 
along the unit normal vector of $\partial D$.  For sufficiently 
smooth $D$, these cases correspond to choosing $X_i$ 
distributed as normally reflecting Brownian motion, 
which is reversible, as in Remark \ref{BMExample}.  For 
these cases, the desired regularity is then standard.  
\end{remark}

\begin{proposition}\label{bndval} 
Suppose that Condition \ref{bwdreg} is satisfied.  With 
reference to Lemma \ref{solest}, assume $U$ is a 
$\{{\cal F}_t^W\}$-adapted, $L^1(\pi )$-valued process and satisfies 
\[U(t,X_i(t))\leq {\Bbb E}[\Gamma_i(t)|W,X_i(t)],\quad t\geq 0.\]
Let $\bar {g}$ be any continuous function on $\bar {D}$ with $\bar {
g}|_{\partial D}=g$.  
Then for any $t>0$,
\[\lim_{\epsilon\to 0}\frac {\int_{\partial_{\epsilon}(D)}|\Phi U
(t,x)-\bar {g}(x)|\pi (dx)}{\pi (\partial_{\epsilon}(D))}=0\]
in $L^1(P)$.
\end{proposition}

\begin{proof}
By Lemma \ref{cndrep} and Jensen's inequality, we have
\begin{eqnarray*}
{\Bbb E}\left[\pi (\partial_{\epsilon}(D))^{-1}\int_{\partial_{\epsilon}
(D)}|\Phi U(t,x)-\bar g(x)|\pi (dx)\right]\leq {\Bbb E}\left[|A_i^
U(t)-\bar g(X_i(t))|\big|X_i(t)\in\partial_{\epsilon}(D)\right].\end{eqnarray*}
Recalling the definition of $A_i^U(t)$, we have
\begin{eqnarray}
|A_i^U(t)-\bar {g}(X_i(t))|&\leq&|g(X_i(\tau_i(t)))-\bar {g}(X_i(
t))|{\bf 1}_{\{\tau_i(t)>0\}}+2\|h\|_{\infty}\vee\|g\|_{\infty}{\bf 1}_{
\{\tau_i(t)=0\}}\non\\
&&\qquad +\int_{\tau_i(t)}^t|G(U(s,X_i(s),X_i(s))A_i^U(s)|ds
\label{agest}\\
&&\qquad +\Big|\int_{D\times (\tau_i(t),t]}\rho (X_i(s),u)W(du,ds
)\Big|.\non\end{eqnarray}
Next, observe that
\begin{eqnarray*}
&&{\Bbb E}\left[|g(X_i(\tau_i(t)))-\bar g(X_i(t))|{\bf 1}_{\{\tau_
i(t)>0\}}\Big|X_i(t)\in\partial_{\epsilon}(D)\right]\\
&&\qquad ={\Bbb E}\left[|g(X_i^{*}(\sigma_i))-\bar g(X_i^{*}(0))|
{\bf 1}_{\{\sigma_i<t\}}\Big|X_i^{*}(0)\in\partial_{\epsilon}(D)\right
]\\
&&\qquad\qquad\leq\sup_{x\in\partial_{\epsilon}K}{\Bbb E}\left[|g
(X_i^{*}(\sigma_i))-\bar g(x)|{\bf 1}_{\{\sigma_i<t\}}\Big|X_i^{*}
(0)=x\right].\end{eqnarray*}
Since $D$ is bounded under Condition \ref{diffcnd}, $\overline D$ is 
compact.  Then there exists $x_0\in\partial D$ and $x_n\to x_0$ so that 
\begin{eqnarray*}
&&\limsup_{\epsilon\to 0}\sup_{x\in\partial_{\epsilon}(D)}{\Bbb E}\left
[|g(X_i^{*}(\sigma_i))-\bar g(x)|{\bf 1}_{\{\sigma_i<t\}}\Big|X_i^{
*}(0)=x\right]\\
&&\qquad\qquad =\lim_{n\to\infty}{\Bbb E}\left[|g(X_i^{*}(\sigma_
i))-\bar g(x_n)|{\bf 1}_{\{\sigma_i<t\}}\Big|X_i^{*}(0)=x_n\right
].\end{eqnarray*}
Continuity of $\bar {g}$ and (\ref{bdyreg}) imply that the limit is 
zero.  A similar argument and (\ref{bdyhit}) show that 
$P\left(\tau_i(t)=0|X_i(t)\in\partial_{\epsilon}(D)\right)$ tends to zero, so that the 
conditional expectations of the first two terms on the 
right of (\ref{agest}) go to zero.  Next, we recall that 
\begin{eqnarray*}
|G(U(s,X_i(s),X_i(s))A_i^U(s)|\leq L_1\left(1+|U(s,X_i(s))|^2\right
)|A_i^U(s)|.\end{eqnarray*}
By the Cauchy-Schwarz inequality,
\begin{eqnarray*}
&&{\Bbb E}\left[\int_{\tau_i(t)}^t|G(U(s,X_i(s),X_i(s))A_i^U(s)|d
s\Big|X_i(t)\in\partial_{\epsilon}(D)\right]\\
&&\leq L_1{\Bbb E}\left[t-\tau_i(t)|X_i(t)\in\partial_{\epsilon}
(D)\right]^{\frac 12}{\Bbb E}\left[\int_0^t\left(1+|U(s,X_i(s))|^
2\right)^2|A_i^U(s)|^2ds|X_i(t)\in\partial_{\epsilon}(D)\right]^{\frac 
12}.\end{eqnarray*}
The second factor in this  inequality is 
uniformly bounded by Lemma \ref{condbd}, and
\begin{eqnarray*}
&&{\Bbb E}\left[t-\tau_i(t)|X_i(t)\in\partial_{\epsilon}(D)\right
]^{\frac 12}\\
&&\qquad\qquad ={\Bbb E}\left[\sigma_i\wedge t|X_i^{*}(0)\in\partial_{
\epsilon}(D)\right]^{\frac 12},\end{eqnarray*}
which tends to zero by \ref{bdyhit}.  Notice that $W$ 
remains white noise when conditioned on $X_i$.  
Using It\^o's isometry and the fact that the $L^2$ norm 
dominates the $L^1$ norm, we have 
\begin{eqnarray*}
&&{\Bbb E}\left[\Big|\int_{\mathbb{U}\times (\tau_i(t),r]}\rho (X_
i(s),u)W(du\times ds)\Big||X_i(t)\in\partial_{\epsilon}(D)\right]\\
&&\qquad\qquad\leq {\Bbb E}\left[\int_{\mathbb{U}\times (\tau_i(t
),r]}\rho^2(X_i(s),u)duds|X_i(t)\in\partial_{\epsilon}(D)\right]^{\frac 
12}\\
&&\qquad\qquad\leq K_2{\Bbb E}\left[t-\tau_i(t)|X_i(t)\in\partial_{
\epsilon}(D)\right]^{\frac 12}.\end{eqnarray*}
As above, the last expression tends to zero by (\ref{bdyhit}).
\end{proof}

\setcounter{equation}{0}

\section{Weak equations with boundary terms} 
\label{secteqbc} At least in general,
uniqueness will not hold for the weak-form stochastic partial 
differential equation (\ref{aceq}) using test functions with 
compact support in $D$. Consequently, to 
obtain an equation that has a unique solution, we must 
enlarge the class of test functions.  We do this in two 
different ways, obtaining two formally different weak-form 
stochastic differential equations; however, for both 
equations, 
we give conditions under which the process given by
the particle representation constructed above is the 
unique solution.  The first extension, which gives the 
simplest form for the stochastic partial differential 
equation, is obtained by taking the test functions to be 
$C_0^2(D)$, the space of  twice continuously differentiable 
functions that vanish on the boundary. The second 
extension is obtained by taking the test functions to be 
${\cal D}({\Bbb A})$, the domain of the generator of the semigroup 
corresponding to the reflecting diffusion giving the 
particle locations. In this section, we prove
 Theorems  \ref{spdeunqI} and \ref{spdeunq2I}.

\subsection{Proof of Theorem \ref{spdeunqI}}\label{C02}
In this  subsection, we apply the results of Section 
\ref{sectbndrymeas}\ to enlarge the class of test functions 
in the definition of the equation 
(\ref{aceq}) to all $\varphi$ in $C^2_0(D)$, the class of $C^2$-functions  
that vanish on $\partial D$. More precisely, we consider test 
functions of the form $\varphi (x,s)$ which are twice continuously 
differentiable in $x$, continuously differentiable in $s$, and 
vanish on $\partial D\times [0,\infty )$. To simplify notation, extend $g$ to all 
of $D$ by setting $g(x)=h(x)$ for $x\in D$.  We assume that 
$g$ is continuous on the boundary, but none of the 
calculations below require this extension to be 
continuous.  

\subsubsection{Derivation of (\ref{weakeqbca1})} We begin 
by showing that the process $V$ we have constructed solves 
(\ref{weakeqbca1}).  Define 
\begin{eqnarray*}
Z_i(t)&=&-g(X_i(t))+\int_0^tG(v(s,X_i(s)),X_i(s))A_i(s)ds\\
&&\qquad +\int_0^tb(X_i(s))ds+\int_{{\Bbb U}\times [0,t]}\rho (X_
i(s),u)W(du\times ds).\end{eqnarray*}
Then
\begin{equation}A_i(t)=g(X_i(t))+Z_i(t)-Z_i(\tau_i(t))=Y_i(t)-Z_i
(\tau_i(t)),\label{decomp2a}\end{equation}
and note that $Y_i$ is a semimartingale.

Since $A_i$, or more precisely, $Z_i\circ\tau_i$, does not appear to 
be a semimartingale, we derive a version of It\^o's 
formula for $A_i\varphi\circ X_i$ from scratch.  Specifically, we 
consider the limit of the telescoping sum 
\begin{eqnarray}
\varphi (X_i(t),t)A_i(t)&=&\varphi (X_i(0),0)A_i(0)
\label{itoderiva}\\
&&\qquad +\sum_k(\varphi (X_i(t_{k+1}),t_{k+1})A_i(t_{k+1})-\varphi 
(X_i(t_k),t_k)A_i(t_k))\nonumber\\
&=&\varphi (X_i(0),0)A_i(0)+\sum A_i(t_k)(\varphi (X_i(t_{k+1}),t_{
k+1})-\varphi (X_i(t_k),t_k))\nonumber\\
&&\qquad +\sum\varphi (X_i(t_{k+1}),t_{k+1})(A_i(t_{k+1})-A_i(t_k
))\nonumber\end{eqnarray}
as the mesh size of the partition $\{t_k,0\leq k\leq m\}$ goes 
to zero.  Since $X_i$ is a semimartingale, $\varphi\circ X_i$ is a semimartingale 
and the second term on the right converges to the usual 
semimartingale integral.  Since (\ref{itoderiva}) is an 
identity, the second sum must also converge.  To 
distinguish
limits of expressions of this form from the usual 
semimartingale integral, we will write
\[\int_0^t\varphi (X_i(s),s)d^{+}A_i(s)=\lim_{\max(t_{k+1}-t_k)\rightarrow 
0}\sum\varphi (X_i(t_{k+1}),t_{k+1})(A_i(t_{k+1})-A_i(t_k)),\]
where the $\{t_k\}$ are partitions of $[0,t]$.  Note that this 
integral differs from the usual semimartingale integral 
in that we evaluate the integrand at the right end point 
of the interval $(t_k,t_{k+1}]$ rather than the left.  Of course 
the $d^{+}$-integral is still bilinear,  so we have
\[\int_0^t\varphi (X_i(s),s)d^{+}A_i(s)=\int_0^t\varphi (X_i(s),s
)d^{+}Y_i(s)-\int_0^t\varphi (X_i(s),s)d^{+}Z_i(\tau_i(s)).\]

Observing that the 
covariation of $\varphi\circ X_i$ and $Y_i$ is zero, applying 
semimartingale integral results, we get limits for 
everything in (\ref{itoderiva}) but 
\[-\sum\varphi (X_i(t_{k+1}),t_{k+1})(Z_i(\tau_i(t_{k+1}))-Z_i(\tau_
i(t_k))).\]
Note that the summands are zero unless $\tau_i(t_{k+1})>\tau_i(t_
k)$ 
which means $X_i(t)\in\partial D$ for some $t\in [t_k,t_{k+1}]$. Breaking this expression into pieces, we have 
\begin{eqnarray*}
&&-\sum\varphi (X_i(t_{k+1}),t_{k+1})(Z_i(\tau_i(t_{k+1}))-Z_i(\tau_
i(t_k)))\\
&&\qquad =\sum\varphi (X_i(t_{k+1}),t_{k+1})(g(X_i(\tau_i(t_{k+1}
)))-g(X_i(\tau_i(t_k)))\\
&&\qquad\qquad -\sum\varphi (X_i(t_{k+1}),t_{k+1})\int^{\tau_i(t_{
k+1})}_{\tau_i(t_k)}G(v(s,X_i(s)),X_i(s))A_i(s)ds\\
&&\qquad\qquad -\sum\varphi (X_i(t_{k+1}),t_{k+1})\int^{\tau_i(t_{
k+1})}_{\tau_i(t_k)}b(X_i(s))ds\\
&&\qquad\qquad -\sum\varphi (X_i(t_{k+1}),t_{k+1})\int_{{\Bbb U}\times 
(\tau_i(t_k),\tau_i(t_{k+1})]}\rho (X_i(s),u)W(du\times ds),\end{eqnarray*}
and it is clear that the last three sums on the right converge to 
zero.  It is not immediately clear that the first sum converges to zero, but it does.

\begin{lemma}
The integral $\int_0^t\varphi (X_i(s),s)d^{+}g(X_i(\tau (s)))=0$ for all $
t\geq 0$. 
\end{lemma}

\begin{proof}
Let $\gamma_i(t)=\inf\{s\geq t:X_i(s)\in\partial D\}$. ``Summing by parts,'' 
\begin{eqnarray*}
&&\sum_{k=0}^{m-1}\varphi (X_i(t_{k+1}),t_{k+1})(g(X_i(\tau_i(t_{
k+1})))-g(X_i(\tau_i(t_k)))\\
&&\quad =-\sum_{k=0}^{m-1}(\varphi (X_i(t_{k+1}),t_{k+1})-\varphi 
(X_i(t_k),t_k))g(X_i(\tau_i(t_k))\\
&&\qquad\qquad +\varphi (X_i(t),t)g(X_i(\tau_i(t))-\varphi (X(0),
0)g(X(0))\\
&&\quad\rightarrow -\int_0^tg(X_i(\tau_i(s-)))d\varphi (X_i(s),s)\\
&&\qquad\qquad +\varphi (X_i(t),t)g(X_i(\tau_i(t)))-\varphi (X_i(
0),0)g(X_i(0))\\
&&\quad =-\int_{\gamma_i(0)}^tg(X_i(\tau_i(s-)))d\varphi (X_i(s),
s)+\varphi (X_i(t),t)g(X_i(\tau_i(t))),\end{eqnarray*}
where the last equality follows from the fact that 
$\tau_i(s)=0$ for $s<\gamma_i(0)$ and $\varphi (X_i(\gamma_i(0)),
0)=0$.
 For each $n$, define
\[U_n(s)=\sum_{k=0}^{\infty}g(X_i(\gamma_i(\frac kn))){\bf 1}_{[\gamma_
i(\frac kn),\gamma_i(\frac {k+1}n))}(s).\]
Since $\varphi (X_i(\gamma_i(t),\gamma_i(t))=0$, 
\begin{eqnarray*}
&\int_{\gamma_i(0)}^tg(X_i(\gamma_i(\frac kn))){\bf 1}_{[\gamma_
i(\frac kn),\gamma_i(\frac {k+1}n))}d\varphi (X_i(s),s)\\
&\qquad\qquad=\left\{\begin{array}{cl}
g(\gamma_i(\frac kn))\varphi (X_i(t),t)&\quad\gamma_i(\frac kn)\leq 
t<\gamma_i(\frac {k+1}n)\\
0&\quad\mbox{\rm \ otherwise}\end{array}
\right.,
\end{eqnarray*}
and we have 
\[\int_{\gamma_i(0)}^tU_n(s-)d\varphi (X_i(s),s)=U_n(t)\varphi (X_
i(t),t).\]
If $\gamma_i(\frac kn)<\gamma_i(\frac {k+1}n)$, then $\frac kn\leq
\gamma_i(\frac kn)<\frac {k+1}n$, and if 
$\gamma_i(\frac kn)\leq s<\gamma_i(\frac {k+1}n)$, then $\gamma_i
(\frac kn)\leq\tau_i(\frac kn)<\frac {k+1}n$.  Consequently,  
$\lim_{n\rightarrow\infty}U_n(s)=g(X_i(\tau_i(s)))$ and
\[g(X_i(\tau_i(t)))\varphi (X_i(t),t)=\lim_{n\rightarrow\infty}\int_{
\gamma_i(0)}^tU_n(s-)d\varphi (X_i(s),s)=\int_{\gamma_i(0)}^tg(\tau_
i(s-))d\varphi (X_i(s),s)\]
proving the lemma.
\end{proof}
Defining 
\[M_{\varphi ,i}(t)=\int_0^t\nabla\varphi (X_i(s),s)^T\sigma (X_i
(s))dB_i(s),\]
we have
\begin{eqnarray*}
\varphi (X_i(t),t)A_i(t)&=&\varphi (X_i(0),0)A_i(0)+\int_0^t\varphi 
(X_i(s),s)d^{+}A_i(s)\\
&&\qquad +\int_0^tA_i(s)dM_{\varphi ,}{}_i(s)+\int_0^tA_i(s)({\Bbb L}
\varphi (X_i(s),s)+\partial \varphi (X_i(s),s))ds\\
&&\qquad +\int_0^tA_i(s-)\nabla\varphi (X_i(s),s)\cdot\eta (X_i(s
))dL_i(s)\\
&=&\varphi (X_i(0),0)g(X_i(0))+\int_0^tA_i(s)\varphi (X_i(s),s)G(
v(s,X_i(s)),X_i(s))ds\\
&&\qquad +\int_0^t\varphi (X_i(s),s)b(X_i(s))ds\\
&&\qquad +\int_{{\Bbb U}\times [0,t]}\varphi (X_i(s),s)\rho (X_i(
s),u)W(du\times ds)\\
&&\qquad +\int_0^tA_i(s)dM_{\varphi ,}{}_i(s)+\int_0^tA_i(s)({\Bbb L}
\varphi (X_i(s),s)+\partial \varphi (X_i(s),s))ds\\
&&\qquad +\int_0^tA_i(s-)\nabla\varphi (X_i(s),s)\cdot\eta (X_i(s
))dL_i(s).\end{eqnarray*}
Applying Theorem \ref{bndryid}, the averaged identity 
becomes (\ref{weakeqbca1}).

\subsubsection{Proof of uniqueness in Theorem 
\ref{spdeunqI}} Fix any $U\in {\cal L}(\pi )$.  We begin by proving 
uniqueness for the equation linearized by replacing $V$ in 
$G$ by $U$, that is, we want to solve
\begin{eqnarray}
\langle\varphi (\cdot ,t),V(t)\rangle&=&\langle\varphi (\cdot ,0)
,h\rangle_{\pi}+\int_0^t\langle\varphi (\cdot ,s)G(U(s,\cdot ),\cdot 
),V(s)\rangle ds\non\\
&&\qquad +\int_0^t\int_{\overline{D}}\varphi (x,s)b(x)\pi (dx)ds
\label{weakeqlina}\\
&&\qquad +\int_{{\Bbb U}\times [0,t]}\int_{\overline{D}}\varphi (x,s)\rho (x,u
)\pi (dx)W(du\times ds)\non\\
&&\qquad +\int_0^t\langle {\Bbb L}\varphi (\cdot ,s)+\partial\varphi 
(\cdot ,s),V(s)\rangle ds\non\\
&&\qquad +\int_0^t\int_{\partial D}g(x)\eta (x)\cdot\nabla\varphi 
(x,s)\beta (dx)ds.\non\end{eqnarray}
We will refer to this equation as the linearized equation 
with input $U$.

Suppose there are two solutions $V_1$ and $V_2$ of (\ref{weakeqlina}), 
and let $\delta V=V_1-V_2$.  Then by linearity,
\begin{eqnarray}
\langle\varphi (\cdot ,t),\delta V(t)\rangle&=&\int_0^t\langle\varphi 
(\cdot ,s)G(U(s,\cdot ),\cdot ),\delta V(s)\rangle ds\label{linid}\\
&&\qquad +\int_0^t\langle {\Bbb L}\varphi (\cdot ,s)+\partial\varphi 
(\cdot ,s),\delta V(s)\rangle ds.\non\end{eqnarray}
The assumption that ${\Bbb A}^0$ is the closure of $\{(\varphi ,{\Bbb L}
\varphi ):\varphi\in {\cal D}_0\}$ 
implies that we can extend this identity to functions $\varphi$ 
which are differentiable in time and satisfy 
$\varphi (\cdot ,s)\in {\cal D}({\Bbb A}^0)$ with ${\Bbb A}^0\varphi$ and $
\partial\varphi$ bounded and continuous.  
We denote by ${\Bbb T}^0$ the semigroup generated by ${\Bbb A}^0$.  Then 
for $\psi\in {\cal D}({\Bbb A}_0)$ and $\varphi (x,s)=r_{\epsilon}
(s){\Bbb T}^0(t-s)\psi (x)$, where 
$0\leq r_{\epsilon}\leq 1$ is continuously differentiable, $r_{\epsilon}
(s)=0$ for $s\geq t$ 
and $r_{\epsilon}(s)=1$, for $s\leq t-\epsilon$, (\ref{linid}) becomes 
\begin{eqnarray}
0=\int_0^t\langle r_{\epsilon}(s){\Bbb T}^0(t-s)\psi G(U(s,\cdot 
),\cdot ),\delta V(s)\rangle ds+\int_0^t\langle\partial r_{\epsilon}(s){\Bbb T}^0(t-s)\psi 
,\delta V(s)\rangle .\label{linid2}\end{eqnarray}
Assuming $r_{\epsilon}(s)\rightarrow {\bf 1}_{[0,t)}(s)$ appropriately, the second term on 
the right of (\ref{linid2}) converges 
to $\langle\psi ,\delta V(t)\rangle$ and hence
\[\langle\psi ,\delta V(t)\rangle =-\int_0^t\langle {\Bbb T}^0(t-
s)\psi (\cdot )G(U(s,\cdot ),\cdot ),\delta V(s)\rangle ds.\]
With reference to Condition \ref{fndcndI},
taking the supremum over $\psi\in {\cal D}({\Bbb A}^0)$ with $|\psi 
|\leq 1$,
\begin{eqnarray*}
\int_{\overline{D}}|\delta v(t,x)|\pi (dx)&\leq&\int_0^t\int_{\overline{D}}|G(U(s,x),x)||\delta 
v(s,x)|\pi (dx)ds\\
&\leq&\int_0^tL_1\int_{\overline{D}}(1+|U(s,x)|^2)|\delta v(s,x)|\pi (dx)ds\\
&\leq&\int_0^tL_1(1+C^2)\int_{\overline{D}}|\delta v(s,x)|\pi (dx)ds\\
&&\qquad +\int_0^tL_1\int_{\overline{D}}{\bf 1}_{\{|U(s,x)|\geq C\}}(1+|U(s,x|^
2)|\delta v(s,x)|\pi (dx)ds.\end{eqnarray*}
Taking expectations of both sides and applying the H\"older 
inequality, we have 
\begin{eqnarray*}
&&\int_{\overline{D}}{\Bbb E}[|\delta v(t,x)|]\pi (dx)\\
&&\qquad\leq\int_0^tL_1(1+C^2)\int_{\overline{D}}{\Bbb E}[|\delta v(s,x)|]\pi 
(dx)ds\\
&&\qquad\qquad +\int_0^tL_1\left(\int_{\overline{D}}P\{|U(s,x)|\geq C\}\pi (dx
)\right)^{1/3}\\
&&\qquad\qquad\qquad\times\left(\int_{\overline{D}}{\Bbb E}[(1+|U(s,x|^2)^3]\pi 
(dx)\right)^{1/3}\left(\int_{\overline{D}}|\delta v(s,x)|^3\pi (dx)\right)^{1/
3}ds\\
&&\qquad\leq\int_0^tL_1(1+C^2)\int_{\overline{D}}{\Bbb E}[|\delta v(s,x)|]\pi 
(dx)ds\\
&&\qquad\qquad +e^{-\varepsilon_TC^2/3}\int_0^tL_1\left(\int_{\overline{D}}{\Bbb E}
[e^{\varepsilon_T|U(s,x)|^2}]\pi (dx)\right)^{1/3}\\
&&\qquad\qquad\qquad\times\left(\int_{\overline{D}}{\Bbb E}[(1+|U(s,x|^2)^3]\pi 
(dx)\right)^{1/3}\left(\int_{\overline{D}}|\delta v(s,x)|^3\pi (dx)\right)^{1/
3}ds,\end{eqnarray*}
for $t\leq T$. As in the proof of Theorem \ref{uniqueAI}, 
Gronwall's inequality implies $\delta v(t,\cdot )=0$ for 
$t\leq T\wedge\frac {\varepsilon_T}{6L_1}$, that is, $t<T$ satisfying 
$\lim_{C\rightarrow\infty}(2L_1(1+C^2)t-\varepsilon_TC^2/3=-\infty$. But local uniqueness 
implies global uniqueness, so we have uniqueness for the 
linearized equation.

To complete the proof, we follow an argument used in 
the proof of Theorem 3.5 of \cite{KX99}.  Let $V$ be the 
solution constructed in Section \ref{sectaddnoise}, and let 
$U$ be another solution of (\ref{weakeqbca1}) in ${\cal L}(\pi )$.  
Taking this $U$ as the input in (\ref{avdef}), 
we can construct a particle 
system with weights $A_i^U$ 
and let $\Phi U$ be given 
by (\ref{Phidef}).  Recall that $\Phi U\in {\cal L}(\pi )$.  Averaging, we 
see that $\Phi U$ satisfies 
\begin{eqnarray}
\langle\varphi (\cdot ,t),\Phi U(t)\rangle&=&\langle\varphi (\cdot 
,0),h\rangle_{\pi}+\int_0^t\langle\varphi (\cdot ,s)G(U(s,\cdot )
,\cdot ),\Phi U(s)\rangle ds\non\\
&&\qquad +\int_0^t\int_{\overline{D}}\varphi (x,s)b(x)\pi (dx)ds\non\\
&&\qquad +\int_{{\Bbb U}\times [0,t]}\int_{\overline{D}}\varphi (x,s)\rho (x,u
)\pi (dx)W(du\times ds)\non\\
&&\qquad +\int_0^t\langle {\Bbb L}\varphi (\cdot ,s)+\partial\varphi 
(\cdot ,s),\Phi U(s)\rangle ds\non\\
&&\qquad +\int_0^t\int_{\partial D}g(x)\eta (x)\cdot\nabla\varphi 
(x,s)\beta (dx)ds,\non\end{eqnarray}
that is, $\Phi U$ satisfies the linearized equation with input 
$U$.  But solutions of the linearized equation are unique, 
and since $U$ is a solution of the nonlinear equation, it is 
also a solution of the linearized equation with input $U$ 
and by uniqueness, we have
 $\Phi U=U$.  But that means $U$ has a particle representation 
which solves 
the same system of equations as the particle 
representation for $V$, and hence Theorem \ref{uniqueAI} 
implies that $U=V$.

\subsection{Proof of Theorem \ref{spdeunq2I}}\label{A}
In this subsection, we take the set of test functions to be ${\cal D}({\Bbb A})$, the 
domain of the generator for the semigroup corresponding 
to the location processes. More precisely, we take
 $\varphi (x,t)$ that are continuously differentiable in $t$ with $
\varphi (x,t)=0$ 
for $t>t_{\varphi}$, and $\varphi (\cdot ,t)\in {\cal D}({\Bbb A}
)$, $t\geq 0$, so that ${\Bbb A}\varphi$ is bounded 
and continuous.  As above, we extend $g$ to all 
of $D$ by setting $g(x)=h(x)$ for $x\in D$.  
Let 
\begin{equation}\gamma_i(s)=\inf\{t>s:X_i(t)\in\partial D\},\label{gamdef}\end{equation}
and note that 
${\bf 1}_{\{\tau_i(t)=0\}}={\bf 1}_{\{\gamma_i(0)>t\}}$. Let $P(dy,ds|x)$ be the conditional distribution of 
$(X_i(\gamma_i(0)),\gamma_i(0))$ given $X_i(0)=x$, and let 
$P\varphi (x)=\int\varphi (y,s)P(dy,ds|x)$. Let $X^{*}$ be the reversed 
process and $\gamma^{*}$ be the first time that $X^{*}$ hits the 
boundary.  

\subsubsection{Derivation of (\ref{wkeqbc1})}
For $\varphi$ satisfying the conditions stated above,
\[M_{\varphi ,i}(t)=\varphi (X_i(t),t)-\int_0^t({\Bbb A}\varphi (
X_i(s),s)+\partial \varphi (X_i(s),s)ds\]
is a $\{{\cal F}_t^{W,X_i}\}$-martingale.  As before, consider
\begin{eqnarray}
\varphi (X_i(t),t)A_i(t)&=&\varphi (X_i(0),0)A_i(0)\label{itoderi
v}\\
&&\qquad +\sum_k(\varphi (X_i(t_{k+1}),t_{k+1})A_i(t_{k+1})-\varphi 
(X_i(t_k),t_k)A_i(t_k))\nonumber\\
&=&\varphi (X_i(0),0)A_i(0)+\sum A_i(t_k)(\varphi (X_i(t_{k+1}),t_{
k+1})-\varphi (X_i(t_k),t_k))\nonumber\\
&&\qquad +\sum\varphi (X_i(t_{k+1}),t_{k+1})(A_i(t_{k+1})-A_i(t_k
))\nonumber\end{eqnarray}
as the mesh size of the partition $\{t_k,0\leq k\leq m\}$ goes 
to zero.  Since $M_{\varphi ,i}$ is a martingale, $\varphi\circ X_
i$ is a semimartingale 
and the second term on the right converges to the usual 
semimartingale integral.  We must evaluate
\[\int_0^t\varphi (X_i(s),s)d^{+}A_i(s),\]
Breaking $A_i$ into its components, the difficulty is again
\[-\sum\varphi (X_i(t_{k+1}),t_{k+1})(Z_i(\tau_i(t_{k+1}))-Z_i(\tau_
i(t_k))),\]
but the limit is different. Breaking this expression into pieces, we have 
\begin{eqnarray*}
&&-\sum\varphi (X_i(t_{k+1}),t_{k+1})(Z_i(\tau_i(t_{k+1}))-Z_i(\tau_
i(t_k)))\\
&&\qquad =\sum\varphi (X_i(t_{k+1}),t_{k+1})(g(X_i(\tau_i(t_{k+1}
)))-g(X_i(\tau_i(t_k)))\\
&&\qquad\qquad -\sum\varphi (X_i(t_{k+1}),t_{k+1})\int^{\tau_i(t_{
k+1})}_{\tau_i(t_k)}G(v(s,X_i(s)),X_i(s))A_i(s)ds\\
&&\qquad\qquad -\sum\varphi (X_i(t_{k+1}),t_{k+1})\int^{\tau_i(t_{
k+1})}_{\tau_i(t_k)}b(X_i(s))ds\\
&&\qquad\qquad -\sum\varphi (X_i(t_{k+1}),t_{k+1})\int_{{\Bbb U}\times 
(\tau_i(t_k),\tau_i(t_{k+1})]}\rho (X_i(s),u)W(du\times ds).\end{eqnarray*}
\setlength{\belowdisplayskip}{.5pc} \setlength{\belowdisplayshortskip}{.5pc}
\setlength{\abovedisplayskip}{.5pc} \setlength{\abovedisplayshortskip}{.5pc}
Summing by parts, we can write the first expression on the right as 
\begin{eqnarray*}
&&\sum_{k=0}^{m-1}\varphi (X_i(t_{k+1}),t_{k+1})(g(X_i(\tau_i(t_{k+1})))-g(X_i(\tau_i(t_k)))\\
&&\quad =-\sum_{k=0}^{m-1}(\varphi (X_i(t_{k+1}),t_{k+1})-\varphi 
(X_i(t_k),t_k))g(X_i(\tau_i(t_k))\\
&&\qquad\qquad +\varphi (X_i(t),t)g(X_i(\tau_i(t))-\varphi (X(0),
0)g(X(0))\\
&&\quad\rightarrow -\int_0^tg(X_i(\tau_i(s-)))dM_{\varphi ,i}(s)-
\int_0^tg(X_i(\tau_i(s)))({\Bbb A}\varphi (X_i(s),s)+\partial \varphi 
(X_i(s),s)ds\\
&&\qquad\qquad +\varphi (X_i(t),t)g(X_i(\tau_i(t)))-\varphi (X_i(
0),0)g(X_i(0)).\end{eqnarray*}%
For the remaining three terms, a summand is nonzero 
only if $\tau_i(t_{k+1})\neq\tau_i(t_k)$ which implies that $X_i$ hits the 
boundary between $t_k$ and $t_{k+1}$. If the mesh size is small,
by continuity, $X_i(t_{k+1})$ must be close to  $X_i(\tau (t_{k+1}
))$. 
Let  ${\cal E}_i(t)$ be the set of complete excursions from the 
boundary in the interval $[0,t]$.  Then the last three 
terms converge to
\begin{eqnarray*}
&&-\sum_{(\alpha ,\beta ]\in {\cal E}_i(t)}\varphi (X_i(\beta ),\beta 
)\int^{\beta}_{\alpha}G(v(s,X_i(s)),X_i(s))A_i(s)ds\\
&&\qquad -\sum_{(\alpha ,\beta ]\in {\cal E}_i(t)}\varphi (X_i(\beta 
),\beta )\int^{\beta}_{\alpha}b(X_i(s))ds\\
&&\qquad -\sum_{(\alpha ,\beta ]\in {\cal E}_i(t)}\varphi (X_i(\beta 
),\beta )\int_{{\Bbb U}\times (\alpha ,\beta ]}\rho (X_i(s),u)W(d
u\times ds).\end{eqnarray*}
Note that for $\gamma_i$ given by (\ref{gamdef}) and
$s\in (\alpha ,\beta )$, $\gamma_i(s)=\beta$, so we can write 
\begin{eqnarray*}
&&\hspace{-2pc}-\int_0^{\infty}\varphi (X_i(s),s)d^{+}Z_i(\tau_i(s))\\
&&=-\int_0^{\infty}g(X_i(\tau_i(s-)))dM_{\varphi ,i}(s)-\int_
0^{\infty}g(X_i(\tau_i(s)))({\Bbb A}\varphi (X_i(s),s)+\partial \varphi 
(X_i(s),s))ds\\
&&\quad\qquad -\varphi (X_i(0),0)g(X_i(0))-\int^{\infty}_0\varphi 
(X_i(\gamma_i(s)),\gamma_i(s))G(v(s,X_i(s)),X_i(s))A_i(s)ds\\
&&\quad\qquad -\int^{\infty}_0\varphi (X_i(\gamma_i(s)),\gamma_i
(s))b(X_i(s))ds\\
&&\quad\qquad -\int_{{\Bbb U}\times (0,\infty )}\varphi (X_i(\gamma_
i(s)),\gamma_i(s))\rho (X_i(s),u)W(du\times ds),\end{eqnarray*}
and we have 
\begin{eqnarray*}
\varphi (X_i(t),t)A_i(t)&=&\varphi (X_i(0),0)A_i(0)+\int_0^t\varphi 
(X_i(s),s)d^{+}A_i(s)\\
&&\quad +\int_0^tA_i(s)dM_{\varphi ,}{}_i(s)+\int_0^t({\Bbb A}\varphi 
(X_i(s),s)+\partial\varphi (X_i(s),s))A_i(s)ds\\
&=&\varphi (X_i(0),0)g(X_i(0))+\int_0^t\varphi (X_i(s),s)G(v(s,X_
i(s)),X_i(s))A_i(s)ds\\
&&\quad +\int_0^t\varphi (X_i(s),s)b(X_i(s))ds\\
&&\quad +\int_{{\Bbb U}\times [0,t]}\varphi (X_i(s),s)\rho (X_i(
s),u)W(du\times ds)\\
&&\quad +\int_0^tA_i(s)dM_{\varphi ,}{}_i(s)+\int_0^t({\Bbb A}\varphi 
(X_i(s),s)+\partial\varphi (X_i(s),s))A_i(s)ds\\
&&\quad -\int_0^t\varphi (X_i(s),s)d^{+}Z_i(\tau_i(s))
\end{eqnarray*}
\vspace{-1pc}
\begin{eqnarray*}
&=&\int_0^t(\varphi (X_i(s),s)-{\bf 1}_{\{\gamma_i(s)\leq t\}}\varphi 
(X_i(\gamma_i(s))),\gamma_i(s)))G(v(s,X_i(s)),X_i(s))A_i(s)ds\\
&&\quad +\int_0^t(\varphi (X_i(s),s)-{\bf 1}_{\{\gamma_i(s)\leq 
t\}}\varphi (X_i(\gamma_i(s))),\gamma_i(s)))b(X_i(s))ds\\
&&\quad +\int_{{\Bbb U}\times [0,t]}(\varphi (X_i(s),s)-{\bf 1}_{
\{\gamma_i(s)\leq t\}}\varphi (X_i(\gamma_i(s))),\gamma_i(s)))\rho 
(X_i(s),u)W(du\times ds)\\
&&\quad +\int_0^t(A_i(s)-g(X_i(\tau_i(s))))dM_{\varphi ,}{}_i(s)\\
&&\quad +\int_0^t({\Bbb A}\varphi (X_i(s),s)+\partial\varphi (X_
i(s),s))(A_i(s)-g(X_i(\tau_i(s))))ds\\
&&\quad +\varphi (X_i(t),t)g(X_i(\tau_i(t)))\end{eqnarray*}
Recalling that $\varphi$ is nonzero only on a finite time 
interval and letting $t=\infty$, averaging gives (\ref{wkeqbc1}).

\subsubsection{Proof of uniqueness in Theorem \ref{spdeunq2I}} 
The proof of uniqueness is essentially the same as the corresponding proof for
Theorem \ref{spdeunqI}.


\renewcommand {\theequation}{A.\arabic{equation}}
\appendix

\setcounter{equation}{0}

\section{Appendix}

\subsection{Gaussian white noise}\label{bcapp1} Let $\mu$ be a 
$\sigma$-finite Borel measure on a complete, separable metric 
space $({\Bbb U},d)$, and let $\ell$ be Lebesgue measure on $[0,\infty 
)$.  
Define ${\cal A}({\Bbb U})=\{A\in {\cal B}({\Bbb U}):\mu (A)<\infty 
\}$ and 
${\cal A}({\Bbb U}\times [0,\infty )=\{C\in {\cal B}({\Bbb U}\times 
[0,\infty )):\mu\times\ell (C)<\infty \}$ with ${\cal A}({\Bbb U}
\times [0,t])$ 
defined similarly.
Then $W=\{W(C):C\in {\cal A}({\Bbb U}\times [0,\infty ))\}$ 
is space-time Gaussian white noise with 
covariance measure $\mu$ if the $W(C)$ are jointly Gaussian, 
with ${\Bbb E}[W(C)]=0$ and ${\Bbb E}[W(C_1)W(C_2)]=\mu\times\ell 
(C_1\cap C_2)$.
In particular, 
for each $A\in {\cal A}({\Bbb U})$, $W(A\times [0,t])$ is Gaussian with 
${\Bbb E}[W(A\times [0,t])]=0$ and $Var(W(A\times [0,t]))=\mu (A)
t$, that is, 
$W(A\times [0,\cdot ])$ is a Brownian motion with parameter $\mu 
(A)$,  
and
\[{\Bbb E}[W(A\times [0,t])W(B\times [0,s])]=\mu (A\cap B)t\wedge 
s.\]

It follows that for disjoint $C_i$ satisfying $\sum_{i=1}^{\infty}
\mu\times\ell (C_i)<\infty$,
\[W(\cup_i^{\infty}C_i)=\sum_{i=1}^{\infty}W(C_i)\]
 almost surely, but the exceptional event of probability 
zero, will in general depend on the sequence $\{C_i\}$, that 
is, while $W$ acts in some ways like a signed measure, it 
is not a random signed measure.

Define the filtration $\{{\cal F}_t^W\}$ by 
\[{\cal F}_t^W=\sigma (W(C):C\in {\cal A}({\Bbb U}\times [0,t]).\]
For $i=1,\ldots ,m$, let $A_i\in {\cal A}({\Bbb U})$ and let $\xi_
i$ be a process 
adapted to $\{{\cal F}_t^W\}$ satisfying ${\Bbb E}[\int_0^t\xi_i(
s)^2ds]<\infty$ for each 
$t>0$.  Define
\[Y(u,t)=\sum_{i=1}^m{\bf 1}_{A_i}(u)\xi_i(t)\]
and
\[Z(t)=\int_{U\times [0,t]}Y(u,s)W(du\times ds)=\sum_{i=1}^m\int_
0^t\xi_i(s)dW(A_i\times [0,s]).\]
Then, from properties of the It\^o integral,
\[{\Bbb E}[Z(t)]=0\qquad [Z]_t=\int_0^t\int_{{\Bbb U}}Y^2(u,s)\mu 
(du)ds\]
\[{\Bbb E}[Z^2(t)]={\Bbb E}[[Z]_t]=\int_0^t\int_{{\Bbb U}}{\Bbb E}
[Y^2(u,s)]\mu (du)ds.\]
By this last identity, the integral can be extended to all 
measurable processes $Y$ that are adapted to $\{{\cal F}_t^W\}$ in the 
sense that $Y$ restricted to $[0,t]\times {\Bbb U}$ is 
${\cal B}([0,t])\times {\cal B}({\Bbb U})\times {\cal F}_t^W$-measurable and satisfy 
\[{\Bbb E}[\int_0^t\int_{{\Bbb U}}Y^2(u,s)\mu (du)ds]<\infty ,\quad 
t>0.\]
By a truncation argument, the integral can be extended 
to adapted $Y$ satisfying 
\[\int_0^t\int_{{\Bbb U}}Y^2(u,s)\mu (du)ds<\infty\quad a.s.,\quad 
t>0.\]
Note that under this last extension, the integral is a 
locally square integrable martingale with 
\[[Z]_t=\int_0^t\int_{{\Bbb U}}Y^2(u,s)\mu (du)ds.\]

\subsection{Measurability of density 
process}\label{measdens}

\begin{lemma}\label{opmeas} 
Let $E$ and ${\Bbb W}$ be complete separable metric spaces, let $
X$ 
be a stationary Markov process in $E$ with no fixed 
points of discontinuity, and let $W$ be a ${\Bbb W}$-valued random 
variable that is independent of $X$.  Let 
$f:[0,\infty )\times D_E[0,\infty )\times {\Bbb W}\rightarrow {\Bbb R}$ be bounded and Borel measurable 
and be nonanticipating in the sense that 
\[f(t,x,w)=f(t,x(\cdot\wedge t-),w),\quad (t,x,w)\in [0,\infty )\times 
D_E[0,\infty )\times {\Bbb W}.\]
Letting $\{{\cal G}_t^X\}$ denote the reverse filtration, 
${\cal G}_t^X=\sigma (X(s-),s\geq t)$, there exists a Borel 
measurable
$g:[0,\infty )\times E\times {\Bbb W}\rightarrow {\Bbb R}$ such that $
\{g(t,X(t-),W),t\geq 0\}$ gives the 
optional projection of $\{f(t,X,W),t\geq 0\}$, that is,
for each reverse stopping 
time $\tau$, ($\{\tau\geq t\}\in\sigma (W)\vee {\cal G}_t^X,t\geq 
0$),
\begin{equation}{\Bbb E}[f(\tau ,X,W)|\sigma (W)\vee {\cal G}_{\tau}^
X]=g(\tau ,X(\tau -),W).\label{eqrev}\end{equation}
\end{lemma}

\begin{proof}
Let $R$ be the collection of bounded, Borel measurable 
functions $f$ for which the assumptions and conclusions 
of the lemma hold.  Then $R$ is closed under bounded, 
pointwise limits of increasing functions and under 
uniform convergence.

For $0\leq t_1<\cdots <t_m$, $f_i\in B([0,\infty )\times E)$, $f_
0\in B({\Bbb W})$, let
\begin{equation}f(t,x,w)=f_0(t,w)\prod f_i(t,x(t_i\wedge t-))\in 
R.\label{h0}\end{equation}
Then letting $\{T^{*}(t)\}$ denote the semigroup for the 
time-reversed process, $g$ can be expressed in terms of 
$\{T^{*}(t)\}$ and the $f_i$.  For example, if $m=2$,
\[g(t,x(t-),w)=f_0(t,w)T^{*}(t-t_2\wedge t)[f_2T^{*}(t_2\wedge t-
t_1\wedge t)f_1](x(t-)).\]

Let $H_0$ be 
the collection of functions of the form (\ref{h0}).  Then 
by the appropriate version of the monotone class 
theorem (for example, Corollary 4.4 in the Appendix of 
\cite{EK86}), $R$ contains all bounded functions that 
are $\sigma (H_0)$ measurable, that is all bounded measurable 
functions such that $f(t,x,w)=f(t,x(\cdot\wedge t-),w)$.
\end{proof}

\begin{corollary}\label{opmon}
For each $T>0$,
\begin{equation}{\Bbb E}[\sup_{0\leq s\leq T}f(s,X,W)|\sigma (W)\vee 
{\cal G}_t^X]\geq g(t,X(t-),W)\quad\forall t\in [0,T]\quad a.s.\label{ineqrev}\end{equation}

\end{corollary}

\begin{proof}
The uniqueness of the optional projection up to 
indistinguishability ensures that if $Z_1(t)\leq Z_2(t)$ for all $t$ 
with probability one, then the optional projection of $Z_1$ 
is less than or equal to the optional projection of $Z_2$ for 
all $t$ with probability one. 
\end{proof}



\begin{thebibliography}{99}

\bibitem{ac72}
Samuel~M. Allen and Joseph~W Cahn.
\newblock Ground state structures in ordered binary alloys with second neighbor
  interactions.
\newblock \emph{Acta Met.}, 20\penalty0 (423), 1972.

\bibitem{ac73}
Samuel~M. Allen and Joseph~W Cahn.
\newblock A correction to the ground state of fcc binary ordered alloys with
  first and second neighbor pairwise interactions.
\newblock \emph{Acta Met.}, 7\penalty0 (1261), 1973.

\bibitem{cat88}
P.~Cattiaux.
\newblock Time reversal of diffusion processes with a boundary condition.
\newblock \emph{Stoch. Proc. Appl}, 28\penalty0 (2):\penalty0 275--292, 1988. \MR{0952834}

\bibitem{CKL14}
Dan Crisan, Thomas~G. Kurtz, and Yoonjung Lee.
\newblock Conditional distributions, exchangeable particle systems, and
  stochastic partial differential equations.
\newblock \emph{Ann. Inst. Henri Poincar\'e Probab. Stat.}, 50\penalty0
  (3):\penalty0 946--974, 2014.
\newblock ISSN 0246-0203.
\newblock URL \url{http://dx.doi.org/10.1214/13-AIHP543}. \MR{3224295}

\bibitem{DI93}
Paul Dupuis and Hitoshi Ishii.
\newblock S{DE}s with oblique reflection on nonsmooth domains.
\newblock \emph{Ann. Probab.}, 21\penalty0 (1):\penalty0 554--580, 1993.
\newblock ISSN 0091-1798.
\newblock URL
  \url{http://links.jstor.org/sici?sici=0091-1798(199301)21:1<554:SWORON>2.0.CO;2-D&origin=MSN}. \MR{1207237}

\bibitem{EK86}
Stewart~N. Ethier and Thomas~G. Kurtz.
\newblock \emph{Markov {P}rocesses: Characterization and {C}onvergence}.
\newblock Wiley Series in Probability and Mathematical Statistics: Probability
  and Mathematical Statistics. John Wiley \& Sons Inc., New York, 1986.
\newblock ISBN 0-471-08186-8. \MR{0838085}

\bibitem{Fri75}
Avner Friedman.
\newblock \emph{Stochastic differential equations and applications. {V}ol. 1}.
\newblock Academic Press [Harcourt Brace Jovanovich, Publishers], New
  York-London, 1975.
\newblock Probability and Mathematical Statistics, Vol. 28. \MR{0494490}

\bibitem{Gra92}
Carl Graham.
\newblock Mc{K}ean-{V}lasov {I}t\^o-{S}korohod equations, and nonlinear
  diffusions with discrete jump sets.
\newblock \emph{Stochastic Process. Appl.}, 40\penalty0 (1):\penalty0 69--82,
  1992.
\newblock ISSN 0304-4149. \MR{1145460}

\bibitem{KR14}
Weining Kang and Kavita Ramanan.
\newblock Characterization of stationary distributions of reflected diffusions.
\newblock \emph{Ann. Appl. Probab.}, 24\penalty0 (4):\penalty0 1329--1374,
  2014. \MR{3210998}

\bibitem{Kol10}
Vassili~N. Kolokoltsov.
\newblock \emph{Nonlinear {M}arkov processes and kinetic equations}, volume 182
  of \emph{Cambridge Tracts in Mathematics}.
\newblock Cambridge University Press, Cambridge, 2010.
\newblock ISBN 978-0-521-11184-3.
\newblock URL \url{http://dx.doi.org/10.1017/CBO9780511760303}. \MR{2680971}

\bibitem{KK10}
Peter~M. Kotelenez and Thomas~G. Kurtz.
\newblock Macroscopic limits for stochastic partial differential equations of
  {M}c{K}ean-{V}lasov type.
\newblock \emph{Probab. Theory Related Fields}, 146\penalty0 (1-2):\penalty0
  189--222, 2010.
\newblock ISSN 0178-8051.
\newblock URL \url{http://dx.doi.org/10.1007/s00440-008-0188-0}. \MR{2550362}

\bibitem{Kur14}
Thomas~G. Kurtz.
\newblock Weak and strong solutions of general stochastic models.
\newblock \emph{Electron. Commun. Probab.}, 19:\penalty0 no. 58, 16, 2014.
\newblock ISSN 1083-589X.
\newblock URL
  \url{http://dx.doi.org.ezproxy.library.wisc.edu/10.1214/ECP.v19-2833}. \MR{3254737}

\bibitem{KP96b}
Thomas~G. Kurtz and Philip~E. Protter.
\newblock Weak convergence of stochastic integrals and differential equations.
  {II}. {I}nfinite-dimensional case.
\newblock In \emph{Probabilistic {M}odels for {N}onlinear {P}artial
  {D}ifferential {E}quations (Montecatini Terme, 1995)}, volume 1627 of
  \emph{Lecture Notes in Math.}, pages 197--285. Springer, Berlin, 1996. \MR{1431303}

\bibitem{KX99}
Thomas~G. Kurtz and Jie Xiong.
\newblock Particle representations for a class of nonlinear {SPDE}s.
\newblock \emph{Stochastic Process. Appl.}, 83\penalty0 (1):\penalty0 103--126,
  1999.
\newblock ISSN 0304-4149. \MR{1705602}

\bibitem{KX04}
Thomas~G. Kurtz and Jie Xiong.
\newblock A stochastic evolution equation arising from the fluctuations of a
  class of interacting particle systems.
\newblock \emph{Commun. Math. Sci.}, 2\penalty0 (3):\penalty0 325--358, 2004.
\newblock ISSN 1539-6746. \MR{2118848}

\bibitem{McK67}
H.~P. McKean, Jr.
\newblock Propagation of chaos for a class of non-linear parabolic equations.
\newblock In \emph{Stochastic {D}ifferential {E}quations ({L}ecture {S}eries in
  {D}ifferential {E}quations, {S}ession 7, {C}atholic {U}niv., 1967)}, pages
  41--57. Air Force Office Sci. Res., Arlington, Va., 1967. \MR{0233437}

\bibitem{MP10}
Peter M{\"o}rters and Yuval Peres.
\newblock \emph{Brownian motion}.
\newblock Cambridge Series in Statistical and Probabilistic Mathematics.
  Cambridge University Press, Cambridge, 2010.
\newblock ISBN 978-0-521-76018-8.
\newblock URL \url{http://dx.doi.org/10.1017/CBO9780511750489}.
\newblock With an appendix by Oded Schramm and Wendelin Werner. \MR{2604525}

\bibitem{Oel84}
Karl Oelschl{\"a}ger.
\newblock A martingale approach to the law of large numbers for weakly
  interacting stochastic processes.
\newblock \emph{Ann. Probab.}, 12\penalty0 (2):\penalty0 458--479, 1984.
\newblock ISSN 0091-1798.
\newblock URL
  \url{http://links.jstor.org/sici?sici=0091-1798(198405)12:2<458:AMATTL>2.0.CO;2-H&origin=MSN}. \MR{0735849}

\bibitem{pw}
G.~Parisi and Yong~Shi Wu.
\newblock Perturbation theory without gauge fixing.
\newblock \emph{Sci. Sinica}, 24\penalty0 (4):\penalty0 483--496, 1981.
\newblock ISSN 0582-236X. \MR{0626795}

\bibitem{SV71}
Daniel~W. Stroock and S.~R.~S. Varadhan.
\newblock Diffusion processes with boundary conditions.
\newblock \emph{Comm. Pure Appl. Math.}, 24:\penalty0 147--225, 1971.
\newblock ISSN 0010-3640. \MR{0277037}

\end{thebibliography}



\ACKNO{The third author 
thanks Cristina Costantini for very useful discussions of 
the properties of reflecting diffusions.}


\end{document}